\documentclass[12pt, twoside, leqno]{article}

\usepackage{amsmath,amsthm}
\usepackage{amssymb}
\usepackage{amsfonts}
\usepackage{eucal}
\usepackage{amsxtra}
\usepackage{makeidx}
\usepackage{xypic}
\usepackage{mathrsfs}
\usepackage{pdflscape}
\usepackage{tabularx}
\usepackage{xy}
\usepackage{MnSymbol}
\usepackage{enumitem}

\newtheorem{theorem}{Theorem}[section]
\newtheorem{corollary}[theorem]{Corollary}
\newtheorem{lemma}[theorem]{Lemma}
\newtheorem{proposition}[theorem]{Proposition}

\theoremstyle{definition}
\newtheorem{definition}[theorem]{Definition}
\newtheorem{remark}[theorem]{Remark}
\newtheorem{example}[theorem]{Example}

\newtheorem*{xrem}{Remark}

\numberwithin{equation}{section}

\begin{document}

\baselineskip=17pt

\title{Decomposition of idempotent 2-cocycles}

\author{Christos Lamprakis\\
E-mail: chrlambr@hotmail.com
\and 
Theodora Theohari-Apostolidi\\
School of Mathematics,  \\ 
Aristotle University of Thessaloniki,\\
Thessaloniki 54124, Greece\\
E-mail: theohari@math.auth.gr}

\date{}

\maketitle

\renewcommand{\thefootnote}{}

\footnote{2020 \emph{Mathematics Subject Classification}: Primary 16S35; Secondary 11S25.}

\footnote{\emph{Key words and phrases}: Weak crossed product algebras, Weak $2$-cocycles, Ordered monoids, idempotent $2$-cocycles.}

\renewcommand{\thefootnote}{\arabic{footnote}}
\setcounter{footnote}{0}

\begin{abstract}
Let $L$ be a finite Galois field extension of $K$ with Galois group $G$. We decompose any idempotent 2-cocycle $f$ using finite sequences of descending two-sided ideals of the corresponding weak crossed product algebra $A_f$. We specialize the results in case $f$ is the corresponding idempotent 2-cocycle $f_r$ for some semilinear map $r:G\rightarrow \Omega$, where $\Omega$ is a multiplicative monoid with minimum element.
\end{abstract}

\section{Introduction}
Let $L$ be a finite Galois field extension of $K$ with Galois group $G$. A function $f:G\times G\longrightarrow L$ is called a normalized weak $2$-cocycle (of $G$ over $L$) if it satisfies the conditions $f(\sigma,\tau)f(\sigma\tau,\rho)=f^{\sigma}(\tau,\rho)f(\sigma,\tau\rho)$, for $\sigma,\tau,\rho\in G$, and $f(1,\sigma)=f(\sigma,1)=1$, for $\sigma\in G$. Associated to a weak $2$-cocycle $f$ there is a $K$-algebra $A_f$ called the weak crossed product algebra associated to $f$, first introduced in \cite{HLS}. The $K$-algebra $A_f$ is defined as an $L$-vector space $A_f=\sum_{\sigma\in G} Lx_{\sigma}$ having the symbols $x_{\sigma}$, $\sigma\in G$, as an $L$-basis  and multiplication defined by the rules $x_{\sigma}l=l^{\sigma}x_{\sigma}$ and $x_{\sigma}x_{\tau}=f(\sigma,\tau)x_{\sigma\tau}$, for $\sigma,\tau\in G$ and $l\in L$. The inertial group $H_f$ or $H$ if it is clear from the text of $f$ is defined as $H(f)=\{\sigma\in G:f(\sigma,\sigma^{-1})\neq 0\}$. Then $A_f=\sum_{\sigma\in H(f)} Lx_{\sigma}+J_f$, where $J_f$ (or $J$ if $f$ is clear from the text) is the Jacobson radical of $A_f$ and the unique maximal two-sided ideal of $A_f$. With the notion $I\triangleleft A_f$ we always mean that $I$ is a two-sided ideal of $A_f$. D.E. Haile in \cite{Haile} gave the structure of the two-sided ideals of $A_f$. In particular he proved that if $I\neq (0)$ is a two-sided ideal of $A_f$, then $I=\sum Lx_{\sigma}$, where the sum is taken over those $\sigma\in G$ such that $x_{\sigma}\in I$ (\cite{Haile}, Lemma 2.2). Moreover if $f$ is an idempotent $2$-cocycle, i.e. taking only the values $0$ and $1$, then every ideal $I$ of $A_f$ is of the form $I=\sum_{x_{\sigma}\in I} I_{\sigma}$, where $I_{\sigma}$  is the ideal of $A_f$ generated by $x_{\sigma}$ (\cite{Haile}, Proposition 2.4).

We denote by $E^2(G,L)$ and $E^2(G,L;H)$ the set of idempotent 2-cocycles and the set of idempotent 2-cocycles with inertial group $H$  respectively. We set $G^*=G\backslash H$. To avoid trivialities we suppose that $G^*\neq\emptyset$ and so $A_f$ is not a simple algebra. Let $\Omega$ be a multiplicative totally ordered monoid with minimum element $1$ satisfying the relations $x<y\Rightarrow xz<yz$ and $zx<zy$, for all $x,y,z\in\Omega$. We denote by $Sl(G)$ the set of functions $r:G\rightarrow \Omega$ satisfying the relations $r(1)=1$ and $r(\sigma\tau)\leq r(\sigma)r(\tau)$ for every $\sigma,\tau\in G$. Let $M_r=\{\sigma\in G:r(\sigma)=1\}$. It was shown in (\cite{LaTh1}, Proposition 4.2 and Theorem 5.2)  that $M_r$ is a group and the function $f_r:G\times G\rightarrow L$ defined by the rule
	\begin{displaymath}
		f_r(\sigma,\tau)=
			\begin{cases}
				1, & \text{if}\;r(\sigma\tau)=r(\sigma)r(\tau),\\
				0, & \text{if}\;r(\sigma\tau)<r(\sigma)r(\tau),
			\end{cases}
	\end{displaymath}
is an element of $E^2(G,L;M_r)$.

The problem that we are interested in is which elements $f$ of $E^2(G,L)$ afford the relation $f=f_r$, for some $r\in Sl(G)$. In this article as a first step, we decompose $f$ in such a way that its constituents can take the form $f_r$.  The decomposition that we are going to demonstrate is intended to simplify the system of equations implied by Proposition 6.10 of \cite{LaTh1}. The main tool that will be used is a new construction of an idempotent 2-cocycle from an existing one by means of a finite sequence of descending two-sided ideals of $A_f$ (Definition \ref{defcocidealser}). We investigate the particular case of the family of ideals $\{J_f,I\}$, where $J_f$ (or just $J$) is the radical of $A_f$, and examine the corresponding algebra. We specialize the results for the case $f=f_r$, for some $r\in Sl(G)$.

We recall some results from \cite{LaTh1}. Let $t+1$ be the nilpotency of $J$. To every $f$ we correspond a partition $\{N_i\}_{i=1}^t$ of $G^*$ defined by the rule $N_k(f)=\{\sigma\in G^*:x_{\sigma}\in J^k\backslash J^{k+1}\}$, for $1\leq k\leq t$. It holds that (\cite{LaTh1}, Lemma 3.1)
	\begin{displaymath}
		N_1(f)=\{\sigma\in G^*:f(\sigma_1,\sigma_2)=0\;\text{for all}\;\sigma_1,\sigma_2\in G^*\;\text{with}\;\sigma_1\sigma_2=\sigma\}.
	\end{displaymath}
An ordered set $(\sigma_1,\ldots,\sigma_k)$ of elements of $N_1(f)$ is called a generator of $\sigma\in G^*$ with respect to $f$, if $x_{\sigma}=x_{\sigma_1}\ldots x_{\sigma_k}$. We denote by $g_{\sigma}$ a generator of $\sigma\in G^*$ and by $\Gamma_f$ the set of all generators of all the elements of $G^*$ with respect to $f$. So the elements of $\Gamma_f$ are words with letters from $N_1(f)$. The product of two generators is defined by concatenation. For any map $s':N_1(f)\longrightarrow \Omega$ let $\psi(g)=\prod_{\sigma\in g} s'(\sigma)$ for any generator $g$. In Proposition 6.10 and Theorem 6.12 of \cite{LaTh1} it was shown that if we assign appropriate values to $s'$ (values that, for any element $\sigma$ of $N_1(f)$, force all the values of $\psi(g_{\sigma})$, for the various generators of $\sigma$, to be equal, are invariant for the class $H\sigma H$ and such that $\psi(g_{\sigma})$ is minimum among all the different words for $\sigma$ with letters from $N_1(f)$), then for a natural extension $s$ of $s'$ to all of $G$ we get that $s\in Sl(G)$ and $f=f_s$ .

\section{Idempotent 2-cocycles arising from a finite sequence of descending ideals}

Let $f\in E^2(G,L)$ and $I_{\sigma}$, $\sigma\in G^*$, as in Introduction and $H$ the inertial group of $f$. For $\sigma,\tau\in G$ and $h,h_1,h_2\in H$ it is easy to prove that $f(h_1\sigma,\tau h_2)=f(\sigma,\tau)$ and $f(\sigma h,\tau)=f(\sigma,h \tau)$. Also if $f(\sigma,\tau)=1$ and $\sigma\in G^*$ or $\tau\in G^*$, then $\sigma\tau\notin H$. Moreover a direct consequence of these formulas is that if $x_{\sigma}\in I$, then $x_{h_1\sigma h_2}\in I$, for every $h_1,h_2\in H$. Also if $x_{h_1\sigma h_2}\in I$, for some $h_1,h_2\in H$, then $x_{\sigma}\in I$. Another useful observation that we will use frequently in the sequel is the fact that if $x_{\tau}\in I_{\sigma}$, then there exist $\rho_1,\rho_2\in G$ such that $x_{\tau}=x_{\rho_1}x_{\sigma}x_{\rho_2}$.

For $g_1,g_2\in\Gamma_f$, the relation $g_1<g_2$, if and only if $g_1$ is an ordered part of $g_2$ is a partial ordering with least element the empty word. We call the Hasse diagram with regard s to this ordering the graph of generators of $f$. Let $Ann(J)$ be the ideal of left-right annihilators of $J$. We state the following proposition which will be used frequently.

\begin{proposition}{\label{propidealgen}}
Let $f\in E^2(G,L;H)$ and $\sigma\in G^*$. Then
	\begin{itemize}
		\item[i.] $\{x_{\tau}\in I_{\sigma}:\tau\in G^*\}=\bigcup_{g_{\sigma}}\{x_{\tau}:\exists g_{\tau}\in\Gamma_f\;\text{and}\; h_1,h_2\in H\;\text{with}\; h_1g_{\sigma}h_2\leq g_{\tau}\}$.
		\item[ii.] If $x_{\sigma}\in Ann(J)$, then $x_{h_1\sigma h_2}\in Ann(J)$, for $h_1,h_2\in H$, and also
		\begin{displaymath}
			I_{\sigma}=\sum_{h_1,h_2\in H} Lx_{h_1\sigma h_2}.
		\end{displaymath}
	\end{itemize}
\end{proposition}
\textbf{Proof}: $i.)$ Let $x_{\tau}\in I_{\sigma}$, for $\tau\in G^*$, and $g_{\sigma}\in \Gamma_f$. There exist $\rho_1,\rho_2\in G$ such that $x_{\tau}=x_{\rho_1}x_{\sigma}x_{\rho_2}$. If $\rho_1,\rho_2\in H$, then $\rho_1g_{\sigma}\rho_2=g_{\rho_1\sigma \rho_2}=g_{\tau}\in\Gamma_f$ and the claim is true. Similarly for the other possible cases that is $\rho_1\in H$ and $\rho_2\notin H$, $\rho_1\notin H$ and $\rho_2\in H$, or $\rho_1,\rho_2\in G^*$.

For the opposite direction, suppose that there exist $g_{\sigma},g_{\tau}\in\Gamma_f$ and $h_1,h_2\in H$ with $h_1g_{\sigma}h_2\leq g_{\tau}$. If $g_{\sigma}=(\sigma_1,\ldots,\sigma_l)$ then $(h_1\sigma_1,\sigma_2,\ldots,\sigma_{l-1},\sigma_l h_2)\leq g_{\tau}$. We remark that
	\begin{displaymath}
		g_{\tau}=
			\begin{cases}
				(\tau_1,\ldots,\tau_a,h_1\sigma_1,\sigma_2,\ldots,\sigma_{l-1},\sigma_lh_2,\tau_{a+1},\ldots,\tau_b),\\
				(\tau_1,\ldots,\tau_a,h_1\sigma_1,\sigma_2,\ldots,\sigma_{l-1},\sigma_lh_2),\\
				(h_1\sigma_1,\sigma_2,\ldots,\sigma_{l-1},\sigma_lh_2,\tau_{1},\ldots,\tau_b),\\
				(h_1\sigma_1,\sigma_2,\ldots,\sigma_{l-1},\sigma_lh_2),
			\end{cases}
	\end{displaymath}
where $\sigma_i,\tau_j\in G^*$, $a,b\geq 1$. We examine the first case and similarly we work for the second and third. We set $g_{\rho_1}=(\tau_1,\ldots,\tau_a)$ and $g_{\rho_2}=(\tau_{a+1},\ldots,\tau_b)$. Then $g_{\tau}=g_{\rho_1}g_{h_1\sigma h_2}g_{\rho_2}=g_{\rho_1}h_1g_{\sigma}h_2g_{\rho_2}=g_{\rho_1 h_1}g_{\sigma}g_{h_2\rho_2}$, hence $x_{\tau}\in I_{\sigma}$. For the fourth case $g_{\tau}=h_1g_{\sigma}h_2$ and so $x_{\tau}=x_{h_1}x_{\sigma}x_{h_2}$. Therefore in any case $x_{\tau}\in I_{\sigma}$ and we conclude that in order to calculate the ideal $I_{\sigma}$ we are restricted on the elements of $G^*$ whose generators contain a generator of an element of the set $H\sigma H$.

$ii.$ Let $x_{\sigma}\in Ann(J)$ and $h_1,h_2\in H$. For any $\tau\in G^*$ it holds that $\tau h_1\in G^*$ and $x_{\tau}x_{h_1\sigma h_2}=x_{\tau h_1}x_{\sigma}x_{h_2}=0$. Similarly $x_{h_1\sigma h_2}x_{\tau}=0$ and hence $x_{h_1\sigma h_2}\in Ann(J)$. To prove the equality let $x_{\tau}\in I_{\sigma}$. From the statement $i)$ there exist $g_{\sigma},g_{\tau}\in\Gamma_f$ and $h_1,h_2\in H$ such that $g_{h_1\sigma h_2}=h_1g_{\sigma}h_2\leq g_{\tau}$. Since $x_{h_1\sigma h_2}\in Ann(J)$ and the generators of the annihilators have maximum length, it follows that $g_{h_1\sigma h_2}=g_{\tau}$ and so $\tau=h_1\sigma h_2$. The opposite direction is obvious.$\square$

The following construction is the main tool which we use throughout this article. All ideals are assumed two-sided.

\begin{definition}{\label{defcocidealser}}
Let $f\in E^2(G,L;H)$ and $\mathbf{I}=\{I_i\}_{i=1}^k$, $k\geq 2$, be a finite sequence of descending ideals of $A_f$, that is  $I_1\supseteq\ldots\supseteq I_k$. We set 
	\begin{displaymath}
		f_{\mathbf{I}}(\sigma,\tau)=
			\begin{cases}
				1, & \sigma\in H\;\text{or}\;\tau\in H,\\
				1, & f(\sigma,\tau)=1\;\text{and}\;x_{\sigma},x_{\tau},x_{\sigma\tau}\in I_i\backslash I_{i+1},\;\text{for some}\;i\in\{1,\ldots,k-1\},\\
				0, & \text{elsewere},
			\end{cases}
	\end{displaymath}
for $\sigma,\tau\in G$. Also, for $\sigma\in G^*$ such that $x_{\sigma}\in I_1$, let $s(\sigma)=max\{a\in\mathbb{N}^*:x_{\sigma}\in I_a, 1\leq a\leq k\}$.
\end{definition}

\begin{definition}{\label{defidempcocorder}}
Inside $E^2(G,L)$ we define the relation
	\begin{displaymath}
		f\leq g\Leftrightarrow \{(\sigma,\tau)\in G\times G:f(\sigma,\tau)=1\}\subseteq \{(\sigma,\tau)\in G\times G:g(\sigma,\tau)=1\}
	\end{displaymath}
which is a partial order.
\end{definition}

\begin{remark}{\label{remcocidealserinf}}
Let $f_{\mathbf{I}}(\sigma,\tau)=1$ for some $\sigma,\tau\in G$. From Definition \ref{defcocidealser} it follows that $f(\sigma,\tau)=1$ or $\sigma\in H$ or $\tau\in H$. In any case $f(\sigma,\tau)=1$ and so $f_{\mathbf{I}}\leq f$. Also for $x_{\sigma},x_{\tau}\in I_1$, if $s(\sigma)=k$ or $s(\tau)=k$ or $s(\sigma\tau)=k$, then $f_{\mathbf{I}}(\sigma,\tau)=0$.
\end{remark}

Our aim is to prove that $f_{\mathbf{I}}$ is an idempotent 2-cocycle. We need the following lemma.

\begin{lemma}{\label{lemfidealsernormal}}
Let $\mathbf{I}=\{I_i\}_{i=1}^k$, $k\geq 2$, be a finite sequence of descending ideals of $A_f$ and $f_{\mathbf{I}}$ the function of Definition \ref{defcocidealser}. It holds that
	\begin{enumerate}[label=\upshape(\roman*), leftmargin=*, widest=iv]
		\item $s(h_1\sigma h_2)=s(\sigma)$, for $x_{\sigma}\in I_1$, $h_1,h_2\in H$,
		\item If  $x_{\sigma},x_{\tau}\in I_1$ and $f(\sigma,\tau)=1$, then $s(\sigma)\leq s(\sigma\tau)$ and $s(\tau)\leq s(\sigma\tau)$,
		\item $f_{\mathbf{I}}(h_1\sigma,\tau h_2)=f_{\mathbf{I}}(\sigma,\tau)$, for $\sigma,\tau\in G$ and $h_1,h_2\in H$,
		\item $f_{\mathbf{I}}(\sigma h,\tau)=f_{\mathbf{I}}(\sigma,h\tau)$, for $\sigma,\tau\in G$ and $h\in H$.
	\end{enumerate}
\end{lemma}

\begin{proof} (i) If $s(\sigma)=k$ the proof is immediate. Suppose that $1\leq s(\sigma)\leq k-1$. Since $x_{\sigma}\in I_{s(\sigma)}$, we have $x_{h_1\sigma h_2}\in I_{s(\sigma)}$. If it was $x_{h_1\sigma h_2}\in I_{s(\sigma)+1}$, then $x_{h_1^{-1}}x_{h_1\sigma h_2}x_{h_2^{-1}}=x_{\sigma}\in I_{s(\sigma)+1}$, a contradiction. Hence $s(h_1\sigma h_2)=s(\sigma)$.

(ii) We have that $x_{\sigma\tau}=x_{\sigma}x_{\tau}\in I_{s(\sigma)}x_{\tau}\subseteq I_{s(\sigma)}$ and so $s(\sigma\tau)\geq s(\sigma)$. Similarly $s(\sigma\tau)\geq s(\tau)$.

(iii) Let $\sigma,\tau\in G$. If $\sigma\in H$, then by definition, $f_{\mathbf{I}}(h_1\sigma,\tau h_2)=f_{\mathbf{I}}(\sigma,\tau)=1$. Similarly if $\tau\in H$. Next let $\sigma,\tau\in G^*$. If $x_{\sigma}\notin I_1$, then $x_{h_1\sigma}\notin I_1$ and so $f_{\mathbf{I}}(h_1\sigma,\tau h_2)=f_{\mathbf{I}}(\sigma,\tau)=0$. Similarly if $x_{\tau}\notin I_1$ or $x_{\sigma\tau}\notin I_1$. Next suppose that $x_{\sigma},x_{\tau},x_{\sigma\tau}\in I_1$. From the statement 1, $s(h_1\sigma)=s(\sigma)$, $s(\tau h_2)=s(\tau)$ and $s(h_1\sigma\tau h_2)=s(\sigma\tau)$. If $s(\sigma)=k$ or $s(\tau)=k$ or $s(\sigma\tau)=k$, then $f_{\mathbf{I}}(h_1\sigma,\tau h_2)=f_{\mathbf{I}}(\sigma,\tau)=0$. Suppose that $1\leq s(\sigma),s(\tau),s(\sigma\tau)\leq k-1$. If $f(\sigma,\tau)=1$ and $s(\sigma)=s(\tau)=s(\sigma\tau)$, then $s(h_1\sigma)=s(\tau h_2)=s(h_1\sigma\tau h_2)$. So $f(h_1\sigma,\tau h_2)=1$ and by definition $f_{\mathbf{I}}(\sigma,\tau)=f_{\mathbf{I}}(h_1\sigma,\tau h_2)=1$. If $f(\sigma,\tau)=0$ or $s(\sigma),s(\tau),s(\sigma\tau)$ are not all equal, then $f_{\mathbf{I}}(\sigma,\tau)=f_{\mathbf{I}}(h_1\sigma,\tau h_2)=0$.
			
(iv) The proof is similar with (iii) above with slight modifications. In particular we examine cases for $s(\sigma),s(\tau)$ and $s(\sigma h\tau)$.
\end{proof}

\begin{theorem}{\label{theidealsercoc}}
Let $f\in E^2(G,L;H)$. The function $f_{\mathbf{I}}$, for some finite sequence of descending ideals $\mathbf{I}=\{I_i\}_{i=1}^k$ of $A_f$, of Definition \ref{defcocidealser} is an element of $E^2(G,L;H)$.
\end{theorem}

\begin{proof} Firstly we prove that the 2-cocycle condition holds for $f_{\mathbf{I}}$. Let $\sigma,\tau,\rho\in G$. If $\sigma\in H$, then by definition $f_{\mathbf{I}}(\sigma,\tau)=f_{\mathbf{I}}(\sigma,\tau\rho)=1$ . Also from Lemma \ref{lemfidealsernormal}(3), $f_{\mathbf{I}}(\sigma\tau,\rho)=f_{\mathbf{I}}(\tau,\rho)$ and the 2-cocycle condition is true. Similarly if $\tau\in H$ or $\rho\in H$.

Next suppose that $\sigma,\tau,\rho\in G^*$. If $f(\sigma,\tau)=0$, then $f_{\mathbf{I}}(\sigma,\tau)=0$. Since $f$ is an idempotent 2-cocycle, either $f(\tau,\rho)$ or $f(\sigma,\tau\rho)$ is zero, so either $f_{\mathbf{I}}(\tau,\rho)$ or $f_{\mathbf{I}}(\sigma,\tau\rho)$ is zero and the 2-cocycle condition is true. The same argument applies when $f(\sigma\tau,\rho)=0$ or $f(\sigma,\tau\rho)=0$ or $f(\tau,\rho)=0$.

Suppose that $f(\sigma,\tau)=f(\sigma\tau,\rho)=f(\tau,\rho)=f(\sigma,\tau\rho)=1$. If $x_{\sigma}\notin I_1$ or $x_{\sigma}\in I_k$, then $f_{\mathbf{I}}(\sigma,\tau)=f_{\mathbf{I}}(\sigma,\tau\rho)=0$. Similarly in the following cases: If $x_{\tau}\notin I_1$ or $x_{\tau}\in I_k$, if $x_{\rho}\notin I_1$ or $x_{\rho}\in I_k$, if $x_{\sigma\tau}\notin I_1$ or $x_{\sigma\tau}\in I_k$, if $x_{\tau\rho}\notin I_1$ or $x_{\tau\rho}\in I_k$ and if $x_{\sigma\tau\rho}\notin I_1$ or $x_{\sigma\tau\rho}\in I_k$. Therefore in all the above cases the 2-cocycle condition holds.

Finally we suppose that  $x_{\sigma},x_{\tau},x_{\rho},x_{\sigma\tau},x_{\tau\rho},x_{\sigma\tau\rho}\in I_1\backslash I_k$. We distinguish cases for the natural numbers $s(\sigma),s(\tau),s(\rho)$. If $s(\sigma)<s(\tau)$, then $s(\sigma)<s(\tau\rho)$ [Lemma \ref{lemfidealsernormal}(2)] and so $f_{\mathbf{I}}(\sigma,\tau)=f_{\mathbf{I}}(\sigma,\tau\rho)=0$. Hence the 2-cocycle condition holds for $f_{\mathbf{I}}$. If $s(\rho)<s(\tau)$, then $s(\rho)<s(\sigma\tau)$ and so $f_{\mathbf{I}}(\sigma\tau,\rho)=f_{\mathbf{I}}(\tau,\rho)=0$. In case $s(\tau)\leq s(\sigma)$ and $s(\tau)\leq s(\rho)$ we distinguish the following cases:
	\begin{enumerate}
		\item $s(\tau)<s(\sigma)\leq s(\rho)$. Then $f_{\mathbf{I}}(\sigma,\tau)=f_{\mathbf{I}}(\tau,\rho)=0$. Similarly if $s(\tau)<s(\rho)\leq s(\sigma)$.
		\item $s(\tau)=s(\sigma)<s(\rho)$. Then $f_{\mathbf{I}}(\tau,\rho)=0$. If $s(\tau)=s(\sigma)=s(\sigma\tau)$, then $f_{\mathbf{I}}(\sigma\tau,\rho)=0$. If $s(\tau)=s(\sigma)<s(\sigma\tau)$, then $f_{\mathbf{I}}(\sigma,\tau)=0$. Similarly if $s(\tau)=s(\rho)<s(\sigma)$.
		\item $s(\tau)=s(\sigma)=s(\rho)$. If $s(\sigma)<s(\sigma\tau\rho)$, then $f_{\mathbf{I}}(\sigma,\tau\rho)=f_{\mathbf{I}}(\sigma\tau,\rho)=0$. If $s(\sigma)=s(\sigma\tau\rho)$, then $s(\sigma)\leq s(\sigma\tau)\leq s(\sigma\tau\rho)$ and so $s(\sigma\tau)=s(\sigma\tau\rho)$ and similarly $s(\tau\rho)=s(\sigma\tau\rho)$. Then $f_{\mathbf{I}}(\sigma,\tau)=f_{\mathbf{I}}(\sigma\tau,\rho)=f_{\mathbf{I}}(\tau,\rho)=f_{\mathbf{I}}(\sigma,\tau\rho)=1$.
	\end{enumerate}
Therefore again in all the above cases the 2-cocycle condition holds. For the inertial group we have that if $\sigma\in H(f_{\mathbf{I}})$ then $f_{\mathbf{I}}(\sigma,\sigma^{-1})=1$. From Remark \ref{remcocidealserinf}, $f(\sigma,\sigma^{-1})=1$ and so $\sigma\in H$.
\end{proof}

For any subgroup $H$ of $G$ we set
	\begin{displaymath}
		f_0(\sigma,\tau)=
			\begin{cases}
				1, & \text{if}\;\sigma\in H\;\text{or}\;\tau\in H,\\
				0, & \text{elsewhere}
			\end{cases}
	\end{displaymath}
the Waterhouse idempotent. The following proposition will be used frequently. The appropriate 2-cocycle $f_0$ (and its inertial group) will be clear from the context.

\begin{proposition}{\label{propidealsersquare}}
$f_{\mathbf{I}}=f_0$ if and only if $I_a^2\subseteq I_{a+1}$, for every $a\in\{1,\ldots,k-1\}$.
\end{proposition}

\begin{proof} We set $H=H(f_{\mathbf{I}})=H(f_0)$ and $G^*=G\backslash H$. First, suppose that $I_a^2\subseteq I_{a+1}$, for every $a\in\{1,\ldots k-1\}$, and let $\sigma,\tau\in G^*$. Let $f(\sigma,\tau)=1$. If $x_{\sigma}\notin I_1$ or $x_{\sigma}\in I_k$ or $x_{\tau}\notin I_1$ or $x_{\tau}\in I_k$ or $x_{\sigma\tau}\notin I_1$ or $x_{\sigma\tau}\in I_k$, then $f_{\mathbf{I}}(\sigma,\tau)=0$ in these cases.

Next suppose that $x_{\sigma},x_{\tau},x_{\sigma\tau}\in I_1\backslash I_k$. If all three natural numbers $s(\sigma)$, $s(\tau)$, $s(\sigma\tau)$ are not equal, then by definition $f_{\mathbf{I}}(\sigma,\tau)=0$. If $s(\sigma)=s(\tau)=s(\sigma\tau)=a$, for some $1\leq a\leq k-1$, then $x_{\sigma\tau}=x_{\sigma}x_{\tau}\in I_{s(\sigma)}I_{s(\tau)}=I_a^2$. From the assumption it follows that $x_{\sigma\tau}\in I_{a+1}=I_{s(\sigma\tau)+1}$, a contradiction, which implies that this case is impossible. Therefore in any case $f_{\mathbf{I}}(\sigma,\tau)=0$, for every $\sigma,\tau\in G^*$.

Next suppose that $f_{\mathbf{I}}=f_0$. Let $a\in\{1,\ldots,k-1\}$. If $I_a^2\neq 0$ and $x_{\sigma}\in I_a^2$ then there exist $x_{\sigma_1},x_{\sigma_2}\in I_a$ such that $f(\sigma_1,\sigma_2)=1$ with $\sigma_1\sigma_2=\sigma$. If $x_{\sigma}\notin I_{a+1}$, then $s(\sigma)\leq a$. Since from Lemma \ref{lemfidealsernormal}(ii) we have $a\leq s(\sigma_1)\leq s(\sigma)\leq a$, it follows that $s(\sigma_1)=s(\sigma)=a$ and similarly $s(\sigma_2)=s(\sigma)=a$. But then $f_{\mathbf{I}}(\sigma_1,\sigma_2)=1$, a contradiction to the assumption. So $x_{\sigma}\in I_{a+1}$ an the proof is completed.
\end{proof}

We will need a new operation.

\begin{definition}{\label{defveeoper}}
For $\{f_i\}_{i=1}^k$ a finite family of elements of $E^2(G,L)$ we define the operation
	\begin{displaymath}
		\bigvee_{i=1}^k f_i(\sigma,\tau)=f_1(\sigma,\tau)\vee\ldots\vee f_k(\sigma,\tau)=
			\begin{cases}
				0, & \text{if}\;f_i(\sigma,\tau)=0,\;\forall i=1,\ldots,k,\\
				1, & \text{if}\;\exists i\in\{1,\ldots,k\}\;\text{such that}\; f_i(\sigma.\tau)=1.
			\end{cases}
	\end{displaymath}
\end{definition}

The set $E^2(G,L)$ in general is not closed under the operation $\vee$. We will encounter some instance where the result of the operation is indeed an idempotent 2-cocycle. We remark that inside $E^2(G,L;H)$ it holds that $f\vee f_0=f_0\vee f=f$, for $H$ the inertial group of $f_0$.

\begin{lemma}{\label{lemiserbreak}}
Let $f\in E^2(G,L;H)$ and $\{I_i\}_{i=1}^k$, $k\geq 3$ be a finite sequence of descending ideals of  $A_f$. Then for $a\in\{2,\ldots,k-1\}$ we get $f_{\{I_1,\ldots,I_a,\ldots, I_k\}}=f_{\{I_1,\ldots,I_a\}}\vee f_{\{I_a,\ldots,I_k\}}$.
\end{lemma}

\begin{proof} We set $f'=f_{\{I_1,\ldots,I_a,\ldots, I_k\}}$, $f_1=f_{\{I_1,\ldots,I_a\}}$ and $f_2=f_{\{I_a,\ldots,I_k\}}$. First let $\sigma,\tau\in G^*$ such that $f'(\sigma,\tau)=1$. Then $f(\sigma,\tau)=1$ and $x_{\sigma},x_{\tau},x_{\sigma\tau}\in I_b\backslash I_{b+1}$, for some $1\leq b\leq k-1$. If $1\leq b\leq a-1$, then $f_1(\sigma,\tau)=1$. If $a\leq b\leq k-1$, then $f_2(\sigma,\tau)=1$.

For the opposite direction, let $\sigma,\tau\in G$ such that $f_1(\sigma,\tau)\vee f_2(\sigma,\tau)=1$. If $\sigma,\tau\in G^*$, then $f(\sigma,\tau)=1$ and $x_{\sigma},x_{\tau},x_{\sigma\tau}\in I_i\backslash I_{i+1}$, for some $1\leq i\leq a-1$ or $x_{\sigma},x_{\tau},x_{\sigma\tau}\in I_i\backslash I_{i+1}$, for some $a\leq i\leq k-1$. Therefore in any case $f'(\sigma,\tau)=1$.
\end{proof}

\begin{corollary}{\label{coriserbreakmulti}}
Let $\{I_i\}_{i=1}^k$, $k\geq 2$  be a finite sequence of descending ideals of $A_f$. By repeated application of Lemma \ref{lemiserbreak} we have
	\begin{displaymath}
		f_{\{I_1,\ldots,I_k\}}=\bigvee_{i=1}^{k-1} f_{\{I_i,I_{i+1}\}}.\quad\square
	\end{displaymath}
\end{corollary}

The following proposition establishes the formula for decomposing an idempotent 2-cocycle. We will need the fact that if $x_{\sigma}\in\sum_{i=1}^k I_i$, then $x_{\sigma}\in I_i$ for some $i\in\{1,\ldots,k\}$. Also if $x_{\sigma}\notin\sum_{i=1}^k I_i$, then $x_{\sigma}\notin I_i$ for every $i\in\{1,\ldots,k\}$.

\begin{proposition}{\label{propseqsumprop}}
Let $I$ and $I_1,\ldots,I_k$, $k\geq 2$, be ideals of $A_f$ such that $I_i\subseteq I$, for every $i\in\{1,\ldots,k\}$. The following identities hold:
	\begin{enumerate}[label=\upshape(\roman*), leftmargin=*, widest=ii]
		\item $f_{\{I,\sum_{i=1}^k I_i\}}=\prod_{i=1}^k f_{\{I,I_i\}}$,
		\item $f_{\{I,\cap_{i=1}^k I_i\}}=\bigvee_{i=1}^k f_{\{I,I_i\}}$.
	\end{enumerate}
\end{proposition}

\begin{proof} We set $P=\sum_{i=1}^k I_i\lhd A_f$. Then from Theorem \ref{theidealsercoc} we get $H(f_{\{I,P\}})=H(f_{\{I,I_i\}})=H(f)=H$, for every $i$. We set $G^*=G\backslash H$. Let $\sigma,\tau\in G^*$. If $f(\sigma,\tau)=0$, then from Remark \ref{remcocidealserinf} we have that $f_{\{I,P\}}(\sigma,\tau)=f_{\{I,I_i\}}(\sigma,\tau)=0$, for every $i$. Next suppose that $f(\sigma,\tau)=1$. If $x_{\sigma}\notin I$ or $x_{\tau}\notin I$ or $x_{\sigma\tau}\notin I$, then by definition $f_{\{I,P\}}(\sigma,\tau)=f_{\{I,I_i\}}(\sigma,\tau)=0$, for every $i$. If $x_{\sigma},x_{\tau},x_{\sigma\tau}\in I$, then:
	\begin{itemize}
		\item If $x_{\sigma\tau}\notin P$, then $x_{\sigma},x_{\tau}\notin P$. Also $x_{\sigma\tau}\notin I_i$, for every $i$, and so $x_{\sigma},x_{\tau}\notin I_i$, for every $i$. It follows that$f_{\{I,P\}}(\sigma,\tau)=f_{\{I,I_i\}}(\sigma,\tau)=1$, for every $i$.
		\item If $x_{\sigma\tau}\in P$, then $f_{\{I,P\}}(\sigma,\tau)=0$. Also there exists $a\in\{1,\ldots,k\}$ such that $x_{\sigma\tau}\in I_a$. Then $f_{\{I,I_a\}}(\sigma,\tau)=0$ and so $\prod_{i=1}^k f_{\{I,I_i\}}(\sigma,\tau)=0$.
	\end{itemize}

Therefore we get the first statement. For the second we set $P=\bigcap_{i=1}^k I_i\lhd A_f$. The proof is identical with the first in all cases except when $f(\sigma,\tau)=1$ and $x_{\sigma},x_{\tau},x_{\sigma\tau}\in I$. In this case if $x_{\sigma\tau}\notin P$, then $f_{\{I,P\}}(\sigma,\tau)=f_{\{I,I_a\}}(\sigma,\tau)=1$. If $x_{\sigma\tau}\in P$, then $f_{\{I,P\}}(\sigma,\tau)=f_{\{I,I_i\}}(\sigma,\tau)=0$ for every $i$. Therefore we get the second statement.
\end{proof}

\section{The crossed product algebra $A_{f_I}$}

\subsection{The idempotent 2-cocycle $f_{\{J,I\}}$}

The special case of Theorem \ref{theidealsercoc} for the family of ideals $\{I_1=J,I_2=I\}$ for some ideal $I$ of $A_f$ is of interest so we restate it as a separate proposition.

\begin{proposition}{\label{propfIcocycle}}
Let $f\in E^2(G,L;H)$ and $I\lhd A_f$. Then the function defined by the rule
	\begin{displaymath}
		f_I(\sigma,\tau)=f_{\{J,I\}}(\sigma,\tau)=
			\begin{cases}
				1, & \text{if}\;\sigma\in H\;\text{or}\;\tau\in H,\\
				1, & \text{if}\;f(\sigma,\tau)=1,x_{\sigma\tau}\notin I,\;\sigma\notin H,\;\tau\notin H,\\
				0, & \text{elsewhere},
			\end{cases}
	\end{displaymath}
is an element of $E^2(G,L;H)$.
\end{proposition}

Despite the new notion, $f_0$ will still denote the Waterhouse idempotent and not the idempotent 2-cocycle $f_{\{J,0\}}$ corresponding to the zero ideal (which is equal to $f$). The following proposition calculates explicitly the set $N_1(f_I)$.

\begin{proposition}{\label{propn1quotient}}
If $I\lhd A_f$, then $N_1(f_I)=N_1(f)\cup \{\sigma\in G^*:x_{\sigma}\in I\}$.
\end{proposition}

\begin{proof}
Let $H$ be the inertial group of both $f$ and $f_I$ by Proposition \ref{propfIcocycle}. We set $G^*=G\backslash H$. First we prove that $N_1(f)\subseteq N_1(f_I)$. For this let any $\sigma\in N_1(f)$. Suppose that there exist $\tau,\rho\in G^*$ with $\tau\rho=\sigma$ such that $f_I(\tau,\rho)=1$. From Remark \ref{remcocidealserinf} it follows that $f(\tau,\rho)=1$ which is impossible. So $f_I(\tau,\rho)=0$, for every $\tau,\rho\in G^*$ with $\tau\rho=\sigma$, and so $\sigma\in N_1(f_I)$.

Next we prove that $\{\sigma\in G^*:x_{\sigma}\in I\}\subseteq N_1(f_I)$.  For this let $\sigma\in G^*$ such that $x_{\sigma}\in I$ and $\tau,\rho\in G^*$ with $\tau\rho=\sigma$. If it was $f_I(\tau,\rho)=1$, then $f(\tau,\rho)=1$ and $x_{\tau}x_{\rho}=x_{\sigma}\in I$. So by definition $f_I(\tau,\rho)=0$, a contradiction.

For the reverse, let $\sigma\in N_1(f_I)$. If there exist $\tau,\rho\in G^*$ such that $\tau\rho=\sigma$ and $f(\tau,\rho)=1$, then $x_{\tau}x_{\rho}=x_{\sigma}\in I$ (since otherwise if $x_{\sigma}\notin I$ it would follow that $f_I(\tau,\rho)=1$, a contradiction). If $f(\tau,\rho)=0$, for every $\tau,\rho\in G^*$ such that $\tau\rho=\sigma$, then $\sigma\in N_1(f)$. Therefore we have proved the proposition.
\end{proof}

Let $x_{\sigma}\in Ann(J)$. We call $\sigma\in G^*$ a \textbf{trivial annihilator} of $f$, if $\sigma\in N_1(f)$. Otherwise we say that $\sigma$ is a \textbf{non-trivial annihilator} of $f$. All annihilators correspond to maximal elements in both graphs of $f$. Intuitively, the classes of the trivial annihilators lie above $H$ in both graphs of $f$. We recall that every ordered part of a generator is also a generator (\cite{LaTh1}, Proposition 6.3). If $f(\sigma,\tau)=1$ for some $\sigma,\tau\in G^*$, then $g_{\sigma}g_{\tau}=g_{\sigma\tau}$. For what follows we set $A_{f_I}=\sum_{\sigma\in G} Ly_{\sigma}$, where $\{y_{\sigma}:\sigma\in G\}$ is an $L$-basis of $A_{f_I}$.

\begin{lemma}{\label{lemfIannih}}
If $I\lhd A_f$ and $x_{\rho}\in I$, then $\rho$ is a trivial annihilator of $f_I$.
\end{lemma}

\begin{proof}
From Proposition \ref{propn1quotient} we know that $\rho\in N_1(f_I)$. The only generator of $\rho$ with respect to $f_I$ is $(\rho)$ (see \cite{LaTh1}, Remark 6.4). It remains to prove that $y_{\rho}$ is an annihilator of $J_{f_I}$. Suppose that it isn't. We extend $(\rho)$ to a generator (with respect to $f_I$) of an annihilator, say $g_{\tau}=(\rho_1,\ldots,\rho,\ldots,\rho_k)\in \Gamma_{f_I}$, $k\geq 1$ [or $(\rho,\rho_1,\ldots,\rho_k)$ or $(\rho_1,\ldots,\rho_k,\rho)$]. Since $y_{\rho_1}\ldots y_{\rho}\ldots y_{\rho_k}=y_{\tau}$, from Remark \ref{remcocidealserinf} it follows that $x_{\rho_1}\ldots x_{\rho}\ldots x_{\rho_k}=x_{\tau}$. But $x_{\rho}\in I$, so $x_{\tau}\in I$. From Proposition \ref{propn1quotient} it would follow that $\tau\in N_1(f_I)$ and so the only generator of $\tau$ with respect to $f_I$ would be $(\tau)$, a contradiction. So $y_{\rho}$ is an annihilator of $J_{f_I}$.
\end{proof}

\begin{theorem}{\label{thealgebramorp}}
Let $I\lhd A_f$.
\begin{enumerate}[label=\upshape(\roman*), leftmargin=*, widest=ii]
		\item Let $a:G\rightarrow L$ be the function defined by the rule
			\begin{displaymath}
				a(\sigma)=
					\begin{cases}
						1, & \text{if}\;x_{\sigma}\notin I,\\
						0, & \text{if}\;x_{\sigma}\in I.
					\end{cases}
			\end{displaymath}	 
Then the function $\phi:A_f\rightarrow  A_{f_I}$ defined by the rule $\phi(x_{\sigma})=a(\sigma)y_{\sigma}$ and extended by linearity is a $K$-algebra homomorphism with $ker(\phi)=I$,
		\item The function $\psi:A_{f_I}\longrightarrow A_f/I$ defined by the rule
			\begin{displaymath}
				\psi(\sum_{\sigma\in G} l_{\sigma}y_{\sigma})=\sum_{\sigma\in G} l_{\sigma}\overline{x}_{\sigma}
			\end{displaymath}
where $\overline{x}_{\sigma}=x_{\sigma}+I$, is a $K$-algebra epimorphism with kernel $\sum_{x_{\sigma}\in I} Ly_{\sigma}$,
		\item $A_{f_I}\cong (A_f/I)\oplus I$ as $L$-modules.
	\end{enumerate}
\end{theorem}

\begin{proof} (i) If $I=0$, then $f_I=f$, $a(\sigma)=1$ for $\sigma\in G$ and the statement is true. So suppose that $I\neq 0$. First we prove that for $\sigma,\tau\in G$ it holds that
	\begin{displaymath}
		f(\sigma,\tau)a(\sigma\tau)=a(\sigma)a(\tau)f_I(\sigma,\tau).
	\end{displaymath}
If $\sigma\in H$, then $f(\sigma,\tau)=f_I(\sigma,\tau)=1$ and $a(\sigma)=1$. In this case we notice that $a(\sigma\tau)=a(\tau)$, since both $x_{\sigma\tau}=x_{\sigma}x_{\tau},x_{\tau}$ are elements of $I$ or both are not. So the equality holds. Similarly if $\tau\in H$. Next suppose that $\sigma,\tau\in G^*$. Suppose that $f(\sigma,\tau)=1$. If $x_{\sigma\tau}=x_{\sigma}x_{\tau}\notin I$, then $x_{\sigma},x_{\tau}\notin I$, $a(\sigma)=a(\tau)=a(\sigma\tau)=1$ and the equality is true. If $x_{\sigma\tau}\in I$, then $a(\sigma\tau)=0$, $f_I(\sigma,\tau)=0$ and again the equality is true. But $\phi(x_{\sigma}x_{\tau})=f(\sigma,\tau)a(\sigma\tau)y_{\sigma\tau}$ and $\phi(x_{\sigma})\phi(x_{\tau})=a(\sigma)\sigma(a(\tau))f_I(\sigma,\tau)y_{\sigma\tau}$ so $\phi$ is $K$-algebra homomorphism.

For the kernel of $\phi$ we have that if $x=\sum l_{\sigma}x_{\sigma}\in I$ for $\sigma\in G$ with $l_{\sigma}\neq 0$, for every $\sigma$, then $x_{\sigma}\in I$ and so $a(\sigma)=0$. Hence $\phi(x)=\sum l_{\sigma}a(\sigma)y_{\sigma}=0$ and so $x\in ker(\phi)$. Also if $x=\sum l_{\sigma}x_{\sigma}\in ker(\phi)$ with $l_{\sigma}\neq 0$ for every $\sigma$, then $\phi(x)=\sum l_{\sigma}a(\sigma)y_{\sigma}=0$ from where $a(\sigma)=0$ for every $\sigma$ and so $x_{\sigma}\in I$ for every $\sigma$. It follows that $x\in I$.

\noindent (ii) $\psi$ is obviously an $L$-module homomorphism. Let $y_{\sigma},y_{\tau}\in A_{f_I}$. Then $\psi(y_{\sigma}y_{\tau})=f_I(\sigma,\tau)x_{\sigma\tau}+I$
and	$\psi(y_{\sigma})\psi(y_{\tau})=\overline{x}_{\sigma}\overline{x}_{\tau}=f(\sigma,\tau)x_{\sigma\tau}+I$. If $f_I(\sigma,\tau)=1$, then $f(\sigma,\tau)=1$ and so $\psi(y_{\sigma}y_{\tau})=\psi(y_{\sigma})\psi(y_{\tau})=x_{\sigma\tau}+I$. Suppose next that $f_I(\sigma,\tau)=0$. If $x_{\sigma\tau}\in I$, then $\psi(y_{\sigma}y_{\tau})=\psi(y_{\sigma})\psi(y_{\tau})=0$. If $x_{\sigma\tau}\notin I$, then by the definition of $f_I$ we must have $f(\sigma,\tau)=0$ and so again  $\psi(y_{\sigma}y_{\tau})=\psi(y_{\sigma})\psi(y_{\tau})=0$. Finally,
	\begin{displaymath}
		ker(\psi)=\{\sum_{\sigma\in G} l_{\sigma}y_{\sigma}:\sum_{\sigma\in G} l_{\sigma}\overline{x}_{\sigma}=0\}=\{\sum_{\sigma\in G} l_{\sigma}y_{\sigma}:l_{\sigma}=0\;\text{or}\;x_{\sigma}\in I\}=\sum_{x_{\sigma}\in I} Ly_{\sigma}.
	\end{displaymath}
	
\noindent (iii) Consider the function $i:A_f\longrightarrow A_{f_I}$ defined by the rule
	\begin{displaymath}
		i(\sum_{\sigma\in G} l_{\sigma}x_{\sigma})=\sum_{\sigma\in G} l_{\sigma}y_{\sigma}
	\end{displaymath}
which is obviously an $L$-module isomorphism. Then the sequence
\[
\xymatrixrowsep{0.15in}
\xymatrixcolsep{0.3in}
\xymatrix{ 0\ar@{->}[r]	& I\ar@{->}[r]^i& A_{f_I}\ar@<0.4ex>[r]^{\psi} & A_f/I \ar@<0.4ex>[l]^{\theta} \ar@{->}[r] & 0}
\]
is short exact, $\psi$ of the statement ii.) is surjective and $ker(\psi)=\sum_{x_{\sigma}\in I} Ly_{\sigma}=i(I)$. Let $\theta:A_f/I\rightarrow A_{f_I}$ be the identity homomorphism when $A_f/I$ is viewed as a subalgebra of $A_{f_I}$ through $\phi$ of the statement $i)$. Then, for any $\sigma\in G$, $\psi\circ \theta(\overline{x}_{\sigma})=\psi(a(\sigma)y_{\sigma})=a(\sigma)\overline{x}_{\sigma}=\overline{x}_{\sigma}$. So $\psi\circ\theta$ is the identity map of $A_f/I$, and the above sequence is split. The result follows.
\end{proof}

\begin{remark}{\label{remafembed}}
Let $\varphi\in End_K(A_f)$. From Theorem \ref{thealgebramorp}, by selecting $I=ker(\varphi)$ it follows that $\phi(A_f)\cong A_f/ker(\phi)=A_f/ker(\varphi)\cong\varphi(A_f)$. So $\varphi(A_f)$ can be viewed as a subalgebra of $A_{f_{ker(\varphi)}}$ and we have $A_{f_{ker(\varphi)}}=\varphi(A_f)\oplus ker(\varphi)$.
\end{remark}

\subsection{Ideals of $A_{f_I}$}

The idempotent 2-cocycle of Definition \ref{defcocidealser} involves a finite family of ideals. It turns out that to study it, we can focus only in two types of idempotent 2-cocycles, namely $f_{\{J,I\}}$ and $f_{\{I,0\}}$. The next proposition is a key ingredient to the aforementioned claim. Let $\{I_i\}_{i=1}^k$ be a finite sequence of descending ideals of $A_f$. Then $\{I_i/I_k\}_{i=1}^k$ is a finite sequence of descending ideals of $A_f/I_k$. We set $A_k=A_{f_{I_k}}=\sum_{\sigma\in G} L y_{\sigma}$. Let $\psi_k:A_k\longrightarrow A_f/I_k$ be the epimorphism of Theorem \ref{thealgebramorp}(ii).

\begin{proposition}{\label{propidealserquot}}
Using the above notion, $f_{\{I_1,\ldots,I_k\}}=(f_{I_k})_{\{P_1,\ldots,P_{k-1},P_k=ker(\psi_k)\}}$, where $\{P_i\}_{i=1}^k$ is a finite sequence of descending ideals of $A_k$ such that $\psi_k(P_i)=I_i/I_k$, for $1\leq i\leq k$.
\end{proposition}

\begin{proof} Let $A_f=\sum_{\sigma\in G} L x_{\sigma}$. We set $f_1=f_{\{I_1,\ldots,I_k\}}$, $f_2=(f_{I_k})_{\{P_1,\ldots,P_{k-1},P_k=ker(\psi_k)\}}$ and $H=H(f)=H(f_1)=H(f_{I_k})=H(f_2)$. Let $\sigma,\tau\in G^*$ such that $f_2(\sigma,\tau)=1$. We have that $\sigma\tau\in G^*$. Then $f_{I_k}(\sigma,\tau)=1$ and $y_{\sigma},y_{\tau},y_{\sigma\tau}\in P_i\backslash P_{i+1}$, for some $1\leq i\leq k-1$. So $f(\sigma,\tau)=1$ and $x_{\sigma},x_{\tau},x_{\sigma\tau}\in J\backslash I_k$ and $\psi_k(y_{\sigma})$, $\psi_k(y_{\tau})$, $\psi_k(y_{\sigma\tau})\in (I_i/I_k)\backslash (I_{i+1}/I_k)$, for some $1\leq i\leq k-1$. Then $f(\sigma,\tau)=1$ and $x_{\sigma},x_{\tau},x_{\sigma\tau}\notin I_k$ and $x_{\sigma}+I_k,x_{\tau}+I_k,x_{\sigma\tau}+I_k\in (I_i/I_k)\backslash (I_{i+1}/I_k)$, for some $1\leq i\leq k-1$. Then $f(\sigma,\tau)=1$ and $x_{\sigma},x_{\tau},x_{\sigma\tau}\in I_i\backslash I_{i+1}$, for some $1\leq i\leq k-1$. It follows that $f_1(\sigma,\tau)=1$. For the opposite direction we note that 
$\psi_k^{-1}(x_{\sigma}+I_k)=\{y_{\sigma}+x:x\in ker(\psi_k)\}$ and we follow the arguments backwards.
\end{proof}

\begin{lemma}{\label{lemtrivannihreplace}}
Let $I_1\triangleleft A_f$ and $I_2=\sum I_{\sigma}$, for some $x_{\sigma}\in I_1$ such that $\sigma$ is trivial annihilator of $f$. Then $f_{\{I_1,I_2\}}=f_{\{I_1,0\}}$.
\end{lemma}

\begin{proof}
From the identity $f_{\{I_1,I_2+0\}}=f_{\{I_1,I_2\}}f_{\{I_1,0\}}$ of Proposition \ref{propseqsumprop}(i) it follows that $f_{\{I_1,I_2\}}\leq f_{\{I_1,0\}}$. Next let $\sigma,\tau\in G^*$ such that $f_{\{I_1,0\}}(\sigma,\tau)=1$. Then $f(\sigma,\tau)=1$ and $x_{\sigma},x_{\tau},x_{\sigma\tau}\in I_1$. It could not be $x_{\sigma}\in I_2$, since $x_{\sigma}$ is not an annihilator of $J$. Similarly $x_{\tau}\notin I_2$. Finally, since $\sigma\tau\notin N_1(f)$ (\cite{LaTh1}, Lemma 3.1) it follows that $x_{\sigma\tau}\notin I_2$.  So $f(\sigma,\tau)=1$ and $x_{\sigma},x_{\tau},x_{\sigma\tau}\in I_1\backslash I_2$. From the definition $f_{\{I_1,I_2\}}(\sigma,\tau)=1$ and so $f_{\{I_1,0\}}\leq f_{\{I_1,I_2\}}$.
\end{proof}

\begin{remark}{\label{remtrivannihreplace}}
Using Corollary \ref{coriserbreakmulti} we can decompose any $f_{\mathbf{I}}\in E^2(G,L)$ into idempotent 2-cocycles of the type $f_{\{I',I\}}$ which in turn, from Proposition \ref{propidealserquot}, take the form $(f_{I})_{\{P,ker(\psi)\}}$, where $\psi:A_{f_I}\rightarrow A_f/I$ and $P$ is an ideal of $A_{f_I}$ containing $ker(\psi)$ such that $\psi(P)=I'/I$. From Lemma \ref{lemfIannih} we know that if $x_{\sigma}\in I$, then $\sigma$ is a trivial annihilator of $f_I$. From Theorem \ref{thealgebramorp}(ii)  and Lemma \ref{lemtrivannihreplace} we have $(f_I)_{\{P,ker(\psi)\}}=(f_I)_{\{P,0\}}$. In essence, the only idempotent 2-cocycles that are needed are of the types $f_{\{J,I\}}$ and $f_{\{I,0\}}$.
\end{remark}

Another useful fact derived from Proposition \ref{propidealserquot} is that repeated applications of Proposition \ref{propfIcocycle} can be avoided.

\begin{proposition}{\label{propquotsequence1}} 
Let $f\in E^2(G,L)$, $I\lhd A_f$, $P\lhd A_{f_I}$ and $\psi:A_{f_I}\longrightarrow A_f/I$ be the algebra epimorphism of Theorem \ref{thealgebramorp}(ii). Then $(f_I)_P=f_{I_1}$, where $I_1$ is an ideal of $A_f$ such that $I_1/I=\psi(P)$.
\end{proposition}

\begin{proof}
We note that $J_{f_I}$ is the unique maximal ideal in $A_{f_I}$ and so by the lattice isomorphism theorem we have $\psi(J_{f_I})=J_f/I$ which is the unique maximal ideal in $A_f/I$.

First we prove the claim for the case that $P$ contains the kernel of $\psi$. Since $\psi(P)$ is an ideal of $A_f/I$, there exists an ideal $I_1$ of $A_f$ containing $I$ such that $\psi(P)=I_1/I$. Since $I_1$ is nilpotent, let $k$ be the smallest positive integer such that $I_1^{2^k}\subseteq I$. Suppose that $k\geq 2$. We have that $\{I_1^{2^i}/I\}_{i=1}^{k-1}$ is a finite sequence of descending ideals of $A_f/I$. Let $\{P_i\}_{i=1}^{k-1}$ be a finite sequence of descending ideals of $A_{f_I}$ containing $ker(\psi)$ such that $\psi(P_i)=I_1^{2^i}/I$, for $i=1,\ldots,k-1$. We set $P_0=P$. Since $(I'/I)^2\subseteq I'^2/I$, for any ideal $I'\supsetneqq I$, for every $i\in\{0,\ldots,k-2\}$ we have
	\begin{displaymath}
		\psi(P_i^2)\subseteq \psi(P_i)^2\subseteq (I_1^{2^i}/I)^2\subseteq (I_1^{2^i})^2/I=I_1^{2^{i+1}}/I=\psi(P_{i+1}).
	\end{displaymath}
It follows that $P_i^2\subseteq P_{i+1}$. For $i=k-1$ we have $\psi(P_{k-1}^2)\subseteq (I_1^{2^{k-1}}/I)^2=0$ and so $P_{k-1}^2\subseteq ker(\psi)$. Then from Proposition \ref{propidealserquot} and Proposition \ref{propidealsersquare} we have
	\begin{eqnarray}
		f_{\{J_f,I_1,I_1^2,\ldots,I_1^{2^{k-1}},I\}} &=& (f_I)_{\{J_{f_I},P_0,P_1,\ldots,P_{k-1},ker(\psi)\}}\nonumber\\
			&=& (f_I)_{\{J_{f_I},P_0\}}\vee (f_I)_{\{P_0,P_1\}}\vee\ldots\vee (f_I)_{\{P_{k-1},ker(\psi)\}}\nonumber\\
			&=&(f_I)_{\{J_{f_I},P\}}=(f_I)_P.\nonumber
	\end{eqnarray}
But $f_{\{J_f,I_1,I_1^2,\ldots,I_1^{2^{k-1}},I\}}=f_{\{J_f,I_1\}}\vee f_{\{I_1,I_1^2\}}\vee \ldots\vee f_{\{I_1^{2^{k-1}},I\}}=f_{I_1}$ as claimed. The cases for $k=0$ $[I_1=I$ and $P=ker(\psi)]$ and $k=1$ ($I_1^2\subseteq I$) are handled accordingly by omitting the irrelevant terms.

Now suppose that $P$ does not contain $ker(\psi)$. Then for the ideal $P'=P+ker(\psi)$, taking into consideration that $(f_I)_{\ker(\psi)}=(f_I)_{\{J_{f_I},0\}}=f_I$ as in Remark \ref{remtrivannihreplace}, from the first part of the proof and Proposition \ref{propseqsumprop}(i) we have $f_I)_P=(f_I)_P (f_I)_{ker(\psi)}=(f_I)_{P'}=f_{I_1}$ where $I_1$ is an ideal of $A_f$ such that $I_1/I=\psi(P')=\psi(P)$.
\end{proof}

\section{Decomposition of idempotent  2 - cocycles using ideals}

Proposition \ref{propseqsumprop}(ii) together with the next proposition the proof of which is immediate,  imply that for every $f\in E^2(G,L;H)$ the subset $\{f_I:I\triangleleft A_f\}$ is a monoid with respect to the operation $\vee$ with unit element $f_0=f_{\{J,J\}}$ and zero element $f=f_{\{J,0\}}$.

\begin{proposition}{\label{propfIidprop}}
Let $\{I_i\}_{i=1}^k$ be a finite family of two-sided ideals of $A_f$. If $\bigcap_{i=1}^k I_i=\{0\}$, then we have that
	\begin{displaymath}
		f=\bigvee_{i=1}^k f_{I_i}.\quad\square
	\end{displaymath}
\end{proposition}

Proposition \ref{propfIidprop} leads to the following proposition.

\begin{proposition}{\label{propfIeqf}}
Let $I\triangleleft A_f$. Then $f_I=f$ if and only if $I=0$ or $I=\sum I_{\sigma}$, for some trivial annihilators $\sigma\in G^*$ of $f$.
\end{proposition}

\begin{proof}
First suppose that $f_I=f$. Since $H(f_I)=H(f)$, from the definition it follows that $x_{\sigma\tau}\notin I$ for every $\sigma,\tau\notin H$ such that $f(\sigma,\tau)=1$. Suppose that $I\neq 0$ and let $x_{\rho}\in I$, for $\rho\in G^*$. We notice that if $\rho\notin N_1(f)$, then there exist $\sigma,\tau\in G^*$ such that $\rho=\sigma\tau$ with $f(\sigma,\tau)=1$. But $x_{\rho}=x_{\sigma\tau}\in I_{\rho}\subseteq I$ and so $f_I(\sigma,\tau)=0$, a contradiction to the assumption that $f_I=f$. So $\rho\in N_1(f)$. If $\rho$ is not an annihilator of $f$, then there exists $\tau\in G^*$ such that $f(\rho,\tau)=1$ or $f(\tau,\rho)=1$. Consider the first case (the second is handled similarly). We have that $x_{\rho\tau}=x_{\rho}x_{\tau}\in I$ and so $f_I(\rho,\tau)=0$, again a contradiction to the assumption. It follows that $\rho$ is a trivial annihilator of $f$, for every $x_{\rho}\in I$. Since $I=\sum_{x_{\rho}\in I} I_{\rho}$, we are done. The opposite direction is a direct consequence of Lemma \ref{lemtrivannihreplace} and the definition of $f_I$.
\end{proof}

\begin{remark}{\label{remequivannihH}}
The relation $\sigma\sim \tau\Leftrightarrow \sigma\in H\tau H$ is an equivalence relation on $G$. We denote by $[\sigma]$ the class of $\sigma\in G$. Let $A$ be the subset of $G$ such that for every $\sigma\in A$, $x_{\sigma}$ is an annihilator of $J$. $A$ is not empty. If $\sigma\in A$, then from Proposition \ref{propidealgen} we have that $H\sigma H\subseteq A$. So $A=\bigcupdot_{i=1}^k H\sigma_i H$, where $A'=\{\sigma_1,\ldots,\sigma_k\}$ is a complete set of representatives of the classes of $A$.\quad$\square$\\
\end{remark}

For the next proposition we need the following observation: If $\sigma,\tau\in G^*$ with $f(\sigma,\tau)=1$, then $x_{\sigma}\notin I_{\sigma\tau}$.

\begin{proposition}{\label{propidelauniquenta}}
Let $f\in E^2(G,L)$. For every $\rho\in G^*$ such that $\rho\notin N_1(f)$, there exists an ideal $I$ of $A_f$ such that $[\rho]$ is the unique class of non-trivial annihilators of $f_I$, with respect to the equivalence relation of Remark \ref{remequivannihH}.
\end{proposition}

\begin{proof}
Let $I=\sum_{x_{\rho}\notin I_{\sigma}} I_{\sigma}$. We remark that $x_{\rho}\notin I$. First we prove that $\rho$ is an annihilator of $f_I$. For this, let $\tau\in G^*$. If $f(\rho,\tau)=0$, then $f_I(\rho,\tau)=0$ (Remark \ref{remcocidealserinf}). Suppose that $f(\rho,\tau)=1$. Then $x_{\rho}\notin I_{\rho\tau}$ and so $x_{\rho\tau}\in I$. From the definition of $f_I$ it follows that $f_I(\rho,\tau)=0$. Similarly we prove that $f_I(\tau,\rho)=0$, for every $\tau\in G^*$. Since by the assumption $\rho\notin N_1(f)$ and also $x_{\rho}\notin I$, from Proposition \ref{propn1quotient} it follows that $\rho\notin N_1(f_I)$ and so $\rho$ is a non-trivial annihilator of $f_I$.

Finally we prove that $[\rho]$ is the unique class of non-trivial annihilators of $f_I$. Suppose that there exists $\sigma\in G^*$ such that $\sigma$ is a non-trivial annihilator of $f_I$ with $[\rho]\neq [\sigma]$. If it was $x_{\sigma}\notin I$, then $x_{\rho}\in I_{\sigma}$. So there exist $\sigma_1,\sigma_2\in G$ such that $x_{\sigma_1}x_{\sigma}x_{\sigma_2}=x_{\rho}$. If both $\sigma_1,\sigma_2\in H$, then $[\rho]=[\sigma]$ contrary to the assumption. So at least one from the two is an element of $G^*$, say $\sigma_1$. Then $f(\sigma_1,\sigma)=1$ and $x_{\sigma_1\sigma}\notin I$ (otherwise, $x_{\rho}\in I$, impossible). So $f_I(\sigma_1,\sigma)=1$, a contradiction to the assumption that $\sigma$ is an annihilator of $f_I$. If it was $x_{\sigma}\in I$ then $\sigma$ would be a trivial annihilator of $f_I$, which is again a contradiction and the proof is completed.\end{proof}

Aljouiee in \cite{Al} studied weak crossed product algebras whose graphs have a unique maximal element i.e. have no trivial annihilators and a unique class of non-trivial annihilators. In particular he showed that such an algebra is Frobenius (Theorem 1.6). In the next theorem we give a procedure to decompose any idempotent 2-cocyle to idempotent 2-cocycles having a unique class of non-trivial annihilators.

\begin{theorem}{\label{thefIdecomp}}
Let $f\in E^2(G,L)$, $f\neq f_0$. Then $f$ has a unique class of non-trivial annihilators or there exist ideals $\{I_i\}_{i=1}^k$ such that $f=\bigvee_{i=1}^k f_{I_i}$, $f_0<f_{I_{i}}<f$ and each $f_{I_i}$ has a unique class of non-trivial annihilators.
\end{theorem}

\begin{proof}
Let $A=\{\rho_1,\ldots,\rho_k\}$ for some $k\geq 1$ be a complete set of representatives of the classes of elements of $G^*$ such that $x_{\rho_i}\in J^2$ ($A\neq\emptyset$ since $f\neq f_0$). If $k=1$, then $[\rho_1]$ is the unique class of non-trivial annihilators. So suppose that $k\geq 2$. We set $I_i=\sum_{x_{\rho_i}\notin I_{\sigma}} I_{\sigma}$. From Proposition \ref{propidelauniquenta}, $f_{I_i}$  has $[\rho_i]$ as unique class of non-trivial annihilators, for every $i=1,\ldots,k$. If $\bigcap_{i=1}^k I_i=\{0\}$, then by Proposition \ref{propfIidprop}, $f=\bigvee_{i=1}^k f_{I_i}$. If $\bigcap_{i=1}^k I_i\neq \{0\}$, then let $x_{\tau}\in \bigcap_{i=1}^k I_i$. If $x_{\tau}\in J^2$, then $\tau\in[\rho_j]$, for some $j\in\{1,\ldots,k\}$ and so $\tau=h_1\rho_j h_2$, $h_1,h_2\in H$. Since $x_{\rho_j}\notin I_j$, we have $x_{h_1}x_{\rho_j}x_{h_2}=x_{\tau}\notin I_j$ which contradicts the assumption that  $x_{\tau}\in \bigcap_{i=1}^k I_i$. So $x_{\tau}\notin J^2$ and so $\tau\in N_1(f)$. If $\tau$ is not an annihilator of $f$, then there exists $\rho\in G^*$ such that $f(\tau,\rho)=1$ [or $f(\rho,\tau)=1$]. Then $x_{\tau\rho}\in J^2$, and so $x_{\tau}\notin I_{x_{\tau\rho}}=I_a$ for some $a\in\{1,\ldots,k\}$ (or $x_{\tau}\notin I_{x_{\rho\tau}})$. But this contradicts the choice of $x_{\tau}$. It follows that for every $x_{\tau}\in\bigcap_{i=1}^k I_i$, $\tau$ is a trivial annihilator of $f$. From Proposition \ref{propseqsumprop} and Proposition \ref{propfIeqf} we have $\bigvee_{i=1}^k f_{I_i}=f_{\bigcap_{i=1}^k I_i}=f_{\sum I_{\tau}}=f$.

We prove the inequality $f_0<f_{I_i}<f$, for every $i=1,\ldots,k$. We know that $f_I\leq f$. First suppose that $\rho_i$ is a non-trivial annihilator of $f$. Since the class $[\rho_i]$ is not unique, let $\rho$ be another non-trivial annihilator, such that $\rho\notin [\rho_i]$. Then from Proposition \ref{propidealannih} we get that $x_{\rho_i}\notin I_{\rho}$ and so $x_{\rho}\in I_i$. Since $\rho\notin N_1(f)$ it follows that there exist $\sigma_1,\sigma_2\in G^*$ such that $f(\sigma_1,\sigma_2)=1$ with $\sigma_1\sigma_2=\rho$. But then $f_{I_i}(\sigma_1,\sigma_2)=0$, which proves that $f_{I_i}<f$. If $\rho_i$ is not an annihilator where $\rho_i\in N_k(f)$ for some $k\geq 2$, then there exists $\tau\in G^*$ such that $f(\rho_i,\tau)=1$ [or $f(\tau,\rho_i)=1$]. Then $x_{\rho_i}\notin I_{\rho_i\tau}$ and so $x_{\rho_i\tau}\in I_i$. Hence that $f_{I_i}(\rho_i,\tau)=0$, and so again $f_{I_i}<f$.

Finally, since $x_{\rho_i}\notin I_i$ and $x_{\rho_i}\in J^2$, it follows that $J^2\subsetneq I_i$ and so, from Proposition \ref{propidealsersquare}, we have $f_{I_i}> f_0$.
\end{proof}

\begin{example}{\label{exdecomp}}
\noindent Let $G=\mathbb{Z}/9\mathbb{Z}$ and $r=\{0,1,2,3,4,1,2,3,3\}\in Sl(G)$ with  set of generators
	\begin{eqnarray}
		\Gamma_{f_r} &=& \{\{(\overline{1})\},\{(\overline{1},\overline{1})\},\{(\overline{1},\overline{1},\overline{1})\},\{(\overline{5},\overline{8}),(\overline{8},\overline{5}),(\overline{1},\overline{1},\overline{1},\overline{1})\},\{(\overline{5})\},\nonumber\\
		&& \{(\overline{1},\overline{5}),(\overline{5},\overline{1})\},\{(\overline{1},\overline{1},\overline{5}),(\overline{1},\overline{5},\overline{1}),(\overline{5},\overline{1},\overline{1})\},\{(\overline{8})\}\}.\nonumber
	\end{eqnarray}
The corresponding table of values and graph of $f_r$ are

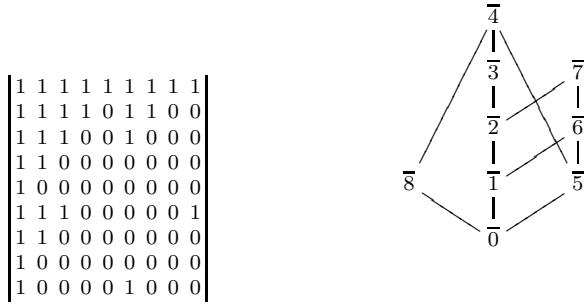
\scriptsize
\[\xymatrixrowsep{0.11in}
\xymatrixcolsep{0.3in}
\setlength{\tabcolsep}{2pt}
\begin{tabular}{m{5cm}m{5cm}}
	\begin{tabular}{|ccccccccc|}
			1&1&1&1&1&1&1&1&1\\
			1&1&1&1&0&1&1&0&0\\
			1&1&1&0&0&1&0&0&0\\
			1&1&0&0&0&0&0&0&0\\
			1&0&0&0&0&0&0&0&0\\
			1&1&1&0&0&0&0&0&1\\
			1&1&0&0&0&0&0&0&0\\
			1&0&0&0&0&0&0&0&0\\
			1&0&0&0&0&1&0&0&0
	\end{tabular} & 
	\xymatrix{ & 	\overline{4}	& \\
		& \overline{3}\ar@{-}[u]& \overline{7}\\
				& \overline{2}\ar@{-}[u]\ar@{-}[ur]& \overline{6}\ar@{-}[u]\\
				\overline{8}\ar@{-}[uuur]&\overline{1}\ar@{-}[u]\ar@{-}[ur]& \overline{5}\ar@{-}[u]\ar@{-}[uuul]\\
					& \overline{0}\ar@{-}[ul]\ar@{-}[u]\ar@{-}[ur] & }
\end{tabular}\]
\normalsize

\noindent The representatives of the classes of the non-trivial annihilators of $f_r$ are $\{\overline{4},\overline{7}\}$. We have that $A=\{\overline{2},\overline{3},\overline{4},\overline{6},\overline{7}\}$. Using the terminology of Theorem \ref{thefIdecomp}, for $\rho_1=\overline{2}$, $I_1=\sum_{x_{\overline{2}}\notin I_{\sigma}} I_{\sigma}=Lx_{\overline{3}}+Lx_{\overline{4}}+Lx_{\overline{5}}+Lx_{\overline{6}}+Lx_{\overline{7}}+Lx_{\overline{8}}$. Similarly for $\rho_2=\overline{3}$, $I_2=Lx_{\overline{4}}+Lx_{\overline{5}}+Lx_{\overline{6}}+Lx_{\overline{7}}+Lx_{\overline{8}}$, for $\rho_3=\overline{4}$, $I_3=Lx_{\overline{6}}+Lx_{\overline{7}}$, for $\rho_4=\overline{6}$, $I_4=Lx_{\overline{2}}+Lx_{\overline{3}}+Lx_{\overline{4}}+Lx_{\overline{7}}+Lx_{\overline{8}}$, and for $\rho_5=\overline{7}$, $I_5=Lx_{\overline{3}}+Lx_{\overline{4}}+Lx_{\overline{8}}$. We know that $f_r=\bigvee_{i=1}^5 f_{I_i}$. Each $f_{I_i}$, $i=1,\ldots,5$, has $[\rho_i]$ as a unique class of non-trivial annihilator. The corresponding graphs are

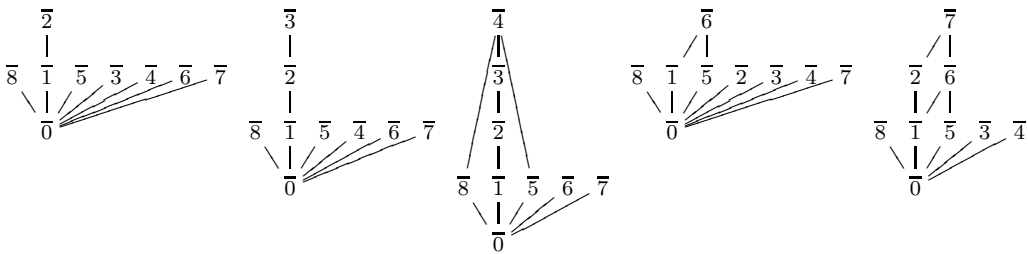
\[\scriptsize
\xymatrixrowsep{0.11in}
\xymatrixcolsep{0.04in}
\xymatrix{	& \overline{2}&  & & & & \\
			\overline{8}&\overline{1}\ar@{-}[u]& \overline{5}& \overline{3} & \overline{4} & \overline{6} & \overline{7}\\
			& \overline{0}\ar@{-}[ul]\ar@{-}[u]\ar@{-}[ur]\ar@{-}[urr]\ar@{-}[urrr]\ar@{-}[urrrr]\ar@{-}[urrrrr] & & & & &}\;
\xymatrix{ 	& \overline{3}&  & & & \\
			& \overline{2}\ar@{-}[u]& & & & \\
			\overline{8}&\overline{1}\ar@{-}[u]& \overline{5}& \overline{4} & \overline{6} & \overline{7} \\
			& \overline{0}\ar@{-}[ul]\ar@{-}[u]\ar@{-}[ur]\ar@{-}[urr]\ar@{-}[urrr]\ar@{-}[urrrr] & & & & }\;
\xymatrix{ & 	\overline{4}	& & & \\
		& \overline{3}\ar@{-}[u]&  & & \\
				& \overline{2}\ar@{-}[u]&  & & \\
				\overline{8}\ar@{-}[uuur]&\overline{1}\ar@{-}[u]& \overline{5}\ar@{-}[uuul] & \overline{6} & \overline{7} \\
					& \overline{0}\ar@{-}[ul]\ar@{-}[u]\ar@{-}[ur]\ar@{-}[urr]\ar@{-}[urrr] & & & }\;
					\xymatrix{	& & \overline{6} & & & & \\
			\overline{8}&\overline{1}\ar@{-}[ur]& \overline{5}\ar@{-}[u]& \overline{2} & \overline{3} & \overline{4} & \overline{7}\\
			& \overline{0}\ar@{-}[ul]\ar@{-}[u]\ar@{-}[ur]\ar@{-}[urr]\ar@{-}[urrr]\ar@{-}[urrrr]\ar@{-}[urrrrr] & & & & &}\;
\xymatrix{	& & \overline{7} & & \\
				& \overline{2}\ar@{-}[ur]& \overline{6}\ar@{-}[u] & & \\
				\overline{8}&\overline{1}\ar@{-}[u]\ar@{-}[ur]& \overline{5}\ar@{-}[u] & \overline{3} & \overline{4}\\
					& \overline{0}\ar@{-}[ul]\ar@{-}[u]\ar@{-}[ur]\ar@{-}[urr]\ar@{-}[urrr] & & & }
\].
\normalsize
\end{example}

\section{Cartesian product of elements of $Sl(G)$}

In this section we specialize the previous results in the case where $f=f_r$, for some $r\in Sl(G)$ taking values in some $\Omega$ as in the introduction. We denote by $Sl(G;H)$ the elements of $Sl(G)$ with $M_r=H$.

\begin{definition}{\label{defidealcartprod}}
Let $r\in Sl(G;H)$ taking values in $\Omega$ and $\mathbf{I}=\{I_i\}_{i=1}^k$, $k\geq 2$, be a finite sequence of descending ideals of $A_{f_r}$. Let $r':G\longrightarrow \times^{k+1} \Omega$ be the function defined by the rule
	\begin{displaymath}
		r_{\mathbf{I}}(\sigma)=
			\begin{cases}
				(\underbrace{r(\sigma),\ldots,r(\sigma)}_{k+1\;\text{times}}), & x_{\sigma}\notin I_1,\\
				(\underbrace{r(\sigma),\ldots,r(\sigma)}_{k-a+1\;\text{times}},\underbrace{1,\ldots,1}_{a\;\text{times}}), & x_{\sigma}\in I_a\backslash I_{a+1}, 1\leq a\leq k-1,\\
				(r(\sigma),\underbrace{1,\ldots,1}_{k\;\text{times}}), & x_{\sigma}\in I_k.\\
			\end{cases}
	\end{displaymath}
\end{definition}

Let $\{\Omega_i,\leq_i\}_{i=1}^k$ be a finite family of multiplicative totally ordered monoids with minimum elements. Then the cartesian product $\times_{i=1}^k \Omega_i=\Omega_1\times\ldots\times\Omega_k$ is a multiplicative monoid with minimum element $1=(1_{\Omega_1},\ldots,1_{\Omega_k})$ totally ordered by the lexicographic relation $(x_1,\ldots,x_k)\leq (y_1,\ldots,y_k)$ if and only if $(x_1,\ldots,x_k)=(y_1,\ldots,y_k)$ or there exists $a\in\{1,\ldots,k\}$ such that, for any $i<a$, $(x_i=y_i$ and $x_a<y_a$). If each $\Omega_i$, $i=1,\ldots,k$ satisfies the relations mentioned in the introduction, then so does $\Omega$.

\begin{theorem}{\label{thercartideal}}
The function of Definition \ref{defidealcartprod} is an element of $Sl(G;H)$.
\end{theorem}

\begin{proof}
Since $x_{h}\notin I_1$ for $h\in H$, we have $r_{\mathbf{I}}(h)=1$. In particular, $r_{\mathbf{I}}(1)=1$. For ease in the calculations we set $I_0=J_{f_r}$ and $I_{k+1}=\{0\}$. Then $r_{\mathbf{I}}$ take the form
	\begin{displaymath}
		r_{\mathbf{I}}(\sigma)=
			\begin{cases}
				(\underbrace{1,\ldots,1}_{k+1\;\text{times}}), & \sigma\in H,\\
				(\underbrace{r(\sigma),\ldots,r(\sigma)}_{k-a+1\;\text{times}},\underbrace{1,\ldots,1}_{a\;\text{times}}), & x_{\sigma}\in I_a\backslash I_{a+1}, 0\leq a\leq k.
			\end{cases}
	\end{displaymath}
First we show that $r_{\mathbf{I}}(h\sigma)=r_{\mathbf{I}}(\sigma h)=r_{\mathbf{I}}(\sigma)$, for $\sigma\in G$ and $h\in H$. If $\sigma\in H$, then $h\sigma,\sigma h\in H$ and so $r_{\mathbf{I}}(h\sigma)=r_{\mathbf{I}}(\sigma h)=r_{\mathbf{I}}(\sigma)=1$. If $x_{\sigma}\in I_a\backslash I_{a+1}$, $0\leq a\leq k$, then also $x_{h\sigma},x_{\sigma h}\in I_a\backslash I_{a+1}$ and so
	\begin{displaymath}
		r_{\mathbf{I}}(h\sigma)=(\underbrace{r(h\sigma),\ldots,r(h\sigma)}_{k-a+1\;\text{times}},\underbrace{1,\ldots,1}_{a\;\text{times}})=r_{\mathbf{I}}(\sigma),
	\end{displaymath}
and similarly $r_{\mathbf{I}}(\sigma h)=r_{\mathbf{I}}(\sigma)$.

Now we show that $r_{\mathbf{I}}(\sigma\tau)\leq r_{\mathbf{I}}(\sigma)r_{\mathbf{I}}(\tau)$, for $\sigma,\tau\in G$. If $\sigma\in H$, then $r_{\mathbf{I}}(\sigma\tau)=r_{\mathbf{I}}(\tau)=r_{\mathbf{I}}(\sigma)r_{\mathbf{I}}(\tau)$ and similarly if $\tau\in H$. If $\sigma\tau\in H$, then $r_{\mathbf{I}}(\sigma\tau)=1\leq r_{\mathbf{I}}(\sigma)r_{\mathbf{I}}(\tau)$. Next let $\sigma,\tau,\sigma\tau\in G^*$. For $\rho\in G^*$ we set $s(\rho)=\{a\in\mathbb{N}:x_{\rho}\in I_a\backslash I_{a+1}, 0\leq a\leq k\}$. We distinguish two cases. First, if $r(\sigma\tau)< r(\sigma)r(\tau)$, for $\sigma,\tau\in G^*$, then
			\begin{displaymath}
				r_{\mathbf{I}}(\sigma\tau)\leq (r(\sigma\tau),\ldots,r(\sigma\tau))< (r(\sigma)r(\tau),1,\ldots,1)\leq r_{\mathbf{I}}(\sigma)r_{\mathbf{I}}(\tau).\nonumber
			\end{displaymath}
Next suppose that $r(\sigma\tau)=r(\sigma)r(\tau)$, for $\sigma,\tau\in G^*$. Then $f_r(\sigma,\tau)=1$. Since $x_{\sigma}\in I_{s(\sigma)}$ it follows that $x_{\sigma}x_{\tau}=x_{\sigma\tau}\in I_{s(\sigma)}$, and so $s(\sigma)\leq s(\sigma\tau)$ and similarly $s(\tau)\leq s(\sigma\tau)$. So we have
		\begin{displaymath}
			 r_{\mathbf{I}}(\sigma\tau)=(\underbrace{r(\sigma\tau),\ldots,r(\sigma\tau)}_{k-s(\sigma\tau)+1},\underbrace{1,\ldots,1}_{s(\sigma\tau)})=(\underbrace{r(\sigma)r(\tau),\ldots,r(\sigma)r(\tau)}_{k-s(\sigma\tau)+1},\underbrace{1,\ldots,1}_{s(\sigma\tau)}).
		\end{displaymath}
Also, for $\sigma,\tau\in G$,
		\begin{displaymath}
			r_{\mathbf{I}}(\sigma)r_{\mathbf{I}}(\tau) = (\underbrace{r(\sigma),\ldots,r(\sigma)}_{k-s(\sigma)+1},\underbrace{1,\ldots,1}_{s(\sigma)})(\underbrace{r(\tau),\ldots,r(\tau)}_{k-s(\tau)+1},\underbrace{1,\ldots,1}_{s(\tau)}).
		\end{displaymath}			
We distinguish two cases, for $\sigma,\tau\in G^*$.
	\begin{itemize}
		\item[a.] $s(\sigma\tau)> max\{s(\sigma),s(\tau)\}$. Then $k-s(\sigma\tau)+1<k-s(\sigma)+1$ and $k-s(\sigma\tau)+1<k-s(\tau)+1$. Since $1<r(\sigma)r(\tau)$, it follows that $r_{\mathbf{I}}(\sigma\tau)<r_{\mathbf{I}}(\sigma)r_{\mathbf{I}}(\tau)$.
		\item[b.] $s(\sigma\tau)=max\{s(\sigma),s(\tau)\}$. If $s(\sigma)=s(\tau)$ then $k-s(\sigma\tau)+1=k-s(\sigma)+1$ and so $r_{\mathbf{I}}(\sigma\tau)=r_{\mathbf{I}}(\sigma)r_{\mathbf{I}}(\tau)$. If $s(\sigma)>s(\tau)$ then $k-s(\sigma\tau)+1=k-s(\sigma)+1$ and $s(\sigma)-s(\tau)>0$. Since $1<r(\tau)$, it follows that $r_{\mathbf{I}}(\sigma\tau)<r_{\mathbf{I}}(\sigma)r_{\mathbf{I}}(\tau)$. Similarly if $s(\sigma)<s(\tau)$.
	\end{itemize}
Finally if $\sigma\in M_{r_{\mathbf{I}}}$, then $r_{\mathbf{I}}(\sigma)=1$. By definition, the only possibility is $r(\sigma)=1$ and so $\sigma\in H$.
\end{proof}

\begin{corollary}{\label{corcareqidealr}}
Let $\mathbf{I}=\{I_i\}_{i=1}^k$, $k\geq 2$, be a finite sequence of descending ideals of $A_{f_r}$. Then $(f_r)_{\mathbf{I}}\leq f_{r_{\mathbf{I}}}\leq f_r$.
\end{corollary}

\begin{proof}
We set $(f_r)_{\mathbf{I}}=f'$. We set $H=H(f')=H(f_r)=M_r=M_{r_{\mathbf{I}}}=H(f_{r_{\mathbf{I}}})$. Let $\sigma,\tau\in G^*$ such that $f'(\sigma,\tau)=1$. Then $f_r(\sigma,\tau)=1$ and $x_{\sigma},x_{\tau},x_{\sigma\tau}\in I_a\backslash I_{a+1}$, for some $1\leq a\leq k-1$. So $r(\sigma\tau)=r(\sigma)r(\tau)$
	\begin{displaymath}
		r_{\mathbf{I}}(\sigma\tau)=(\underbrace{r(\sigma),\ldots,r(\sigma)}_{k-a+1\;\text{times}},\underbrace{1,\ldots,1}_{a\;\text{times}})(\underbrace{r(\tau),\ldots,r(\tau)}_{k-a+1\;\text{times}},\underbrace{1,\ldots,1}_{a\;\text{times}})=r_{\mathbf{I}}(\sigma)r_{\mathbf{I}}(\tau).
	\end{displaymath}
So $f_{r_{\mathbf{I}}}(\sigma,\tau)=1$ which proves the first part of the inequality. For the second part, if $f_r(\sigma,\tau)=0$, then $r(\sigma\tau)<r(\sigma)r(\tau)$ and as in the proof of Theorem \ref{thercartideal}, we deduce that $r_{\mathbf{I}}(\sigma\tau)<r_{\mathbf{I}}(\sigma)r_{\mathbf{I}}(\tau)$ and so $f_{r_{\mathbf{I}}}(\sigma,\tau)=0$. Hence $f_{r_{\mathbf{I}}}\leq f_r$.\end{proof}

\begin{proposition}{\label{propeqcaridealr}}
Let $\mathbf{I}=\{I_i\}_{i=1}^k$, $k\geq 2$, be a finite sequence of descending ideals of $A_{f_r}$ such that $I_1=J_f$ and $I_k=0$. Then $(f_r)_{\mathbf{I}}=f_{r_{\mathbf{I}}}$.
\end{proposition}

\begin{proof}
We set $f_1=(f_r)_{\mathbf{I}}$ and $f_2=f_{r_{\mathbf{I}}}$. We must show that $f_2\leq f_1$. Let $\sigma,\tau\in G$ such that $f_2(\sigma,\tau)=1$. Let $\sigma,\tau\in G^*$. We have that $r_{\mathbf{I}}(\sigma\tau)=r_{\mathbf{I}}(\sigma)r_{\mathbf{I}}(\tau)$. From the proof of Theorem \ref{thercartideal} (case b), we deduce that $r(\sigma\tau)=r(\sigma)r(\tau)$ and $x_{\sigma},x_{\tau},x_{\sigma\tau}\in J_f\backslash I_1$ or $x_{\sigma},x_{\tau},x_{\sigma\tau}\in I_a\backslash I_{a+1}$ for some $1\leq a\leq k-1$ or $x_{\sigma},x_{\tau},x_{\sigma\tau}\in I_k$. Since $I_1=J_f$ and $I_k=0$, the only possibility is $x_{\sigma},x_{\tau},x_{\sigma\tau}\in I_a\backslash I_{a+1}$, for some $1\leq a\leq k-1$. So $f_r(\sigma,\tau)=1$ and by Definition \ref{defcocidealser} $f_1(\sigma,\tau)=1$ as required.
\end{proof}

\begin{theorem}{\label{theridealr}}
Let $r\in Sl(G)$ in any $\Omega$ and $\mathbf{I}=\{I_i\}_{i=1}^k$ be a finite sequence of descending ideals of $A_f$. There exists $r'\in Sl(G)$ in $\times ^l\Omega$, $l\in\mathbb{N}^*$, such that $(f_r)_{\mathbf{I}}=f_{r'}$. 
\end{theorem}

\begin{proof}
Since $J,I_k$ are nilpotent, let $a,b$ be the smallest positive integers such that $J^{2^a}\subseteq I_1$ and $I_k^{2^b}=0$, $a\geq 2$ and $b\geq 1$. From Proposition \ref{propidealsersquare} we notice that $J^2\subseteq J^2+I_1$ and so $(f_r)_{\{J,J^2+I_1\}}=f_0$. Also $(J^{2^{i-1}}+I_1)^2\subseteq (J^{2^{i-1}})^2+I_1=J^{2^i}+I_1$, for $i\in \{2,\ldots,a\}$ and so $(f_r)_{\{J^{2^{i-1}}+I_1,J^{2^i}+I_1\}}=f_0$. Finally, $(I_k^{2^{i-1}})^2\subseteq I_k^{2^i}$, for $i\in \{1,\ldots,b\}$ and so $(f_r)_{\{I_k^{2^{i-1}},I_k^{2^i}\}}=f_0$. Then from Corollary \ref{coriserbreakmulti} we have
\small
	\begin{equation}
		\begin{split}
			  (f_r)&_{\{J,J^2+I_1,\ldots,J^{2^{a-1}}+I_1,I_1,\ldots,I_k,I_k^2,\ldots,I_k^{2^{b-1}},I_k^{2^b}\}}\\
			=(f_r)&_{\{J,J^2+I_1\}}\vee \ldots\vee (f_r)_{\{J^{2^{a-1}}+I_1,I_1\}}\vee (f_r)_{\{I_1,\ldots,I_k\}} \vee (f_r)_{\{I_k,I_k^2\}}\vee\ldots\vee(f_r)_{\{I_k^{2^{b-1}},0\}}\\
			=(f_r)&_{\mathbf{I}}.\nonumber
		\end{split}
	\end{equation}
\normalsize
From Proposition \ref{propeqcaridealr} it follows that the Theorem is true for
	\begin{displaymath}
		r'=r_{\{J,J^2+I_1,\ldots,J^{2^{a-1}}+I_1,I_1,\ldots,I_k,I_k^2,\ldots,I_k^{2^{b-1}},0\}}.
	\end{displaymath}
\noindent There are $a+k+b-2$ ideals so $l=a+k+b-1$. If $a=0$ ($J=I_1$) or $a=1$ ($J^2\subseteq I_1$) or $b=0$ ($I_k=0$) the proof is identical by omitting the irrelevant terms.
\end{proof}

\begin{example}{\label{exdecompr}}
We return to Example \ref{exdecomp}. We notice that $I_3^2=I_5^2=0$, $I_1^2=I_2^2=I_4^2=Lx_{\overline{4}}\neq 0$ and $I_1^4=I_2^4=I_4^4=0$. We set $r_i=r_{\{J,I_i,0\}}$ for $i=3,5$ and $r_i=r_{\{J,I_i,I_i^2,0\}}$ for $i=1,2,4$. Then $f_{I_i}=(f_r)_{\{J,I_i\}}=(f_r)_{\{J,I_i,0\}}=f_{r_i}$ for $i=3,5$ and $f_{I_i}=(f_r)_{\{J,I_i,I_i^2,0\}}=f_{r_i}$ for $i=1,2,4$. More specifically, for $i=1$, we have that $x_{\overline{1}},x_{\overline{2}}\in J\backslash I_1$, $x_{\overline{3}},x_{\overline{5}},x_{\overline{6}},x_{\overline{7}},x_{\overline{8}}\in I_1\backslash I_1^2$, $x_{\overline{4}}\in I_1^2$, for $i=2$, $x_{\overline{1}},x_{\overline{2}},x_{\overline{3}}\in J\backslash I_2$, $x_{\overline{5}},x_{\overline{6}},x_{\overline{7}},x_{\overline{8}}\in I_2\backslash I_2^2$, $x_{\overline{4}}\in I_2^2$, for $i=3$, $x_{\overline{1}},x_{\overline{2}},x_{\overline{3}},x_{\overline{4}},x_{\overline{5}},x_{\overline{8}}\in J\backslash I_3$, $x_{\overline{6}},x_{\overline{7}}\in I_3$, for For $i=4$, we have that $x_{\overline{1}},x_{\overline{5}},x_{\overline{6}}\in J\backslash I_4$, $x_{\overline{2}},x_{\overline{3}},x_{\overline{7}},x_{\overline{8}}\in I_4\backslash I_4^2$, $x_{\overline{4}}\in I_4^2$, and for $i=5$, we have that $x_{\overline{1}},x_{\overline{2}},x_{\overline{5}},x_{\overline{6}},
		x_{\overline{7}}\in J\backslash I_5$, $x_{\overline{3}},x_{\overline{4}},x_{\overline{8}}\in I_5$. So
		
\small
	\begin{displaymath}
		\begin{array}{|c|c|c|c|c|c|c|}
			\hline
			\sigma & r(\sigma) & r_1(\sigma) & r_2(\sigma) & r_3(\sigma) & r_4(\sigma) & r_5(\sigma) \\
			\hline
			\overline{0} & 0 & (0,0,0,0,0) & (0,0,0,0,0) & (0,0,0,0) & (0,0,0,0,0) & (0,0,0,0) \\
			\overline{1} & 1 & (1,1,1,1,0) & (1,1,1,1,0) & (1,1,1,0) & (1,1,1,1,0) & (1,1,1,0) \\
			\overline{2} & 2 & (2,2,2,2,0) & (2,2,2,2,0) & (2,2,2,0) & (2,2,2,0,0) & (2,2,2,0) \\
			\overline{3} & 3 & (3,3,3,0,0) & (3,3,3,3,0) & (3,3,3,0) & (3,3,3,0,0) & (3,3,0,0) \\
			\overline{4} & 4 & (4,4,0,0,0) & (4,4,0,0,0) & (4,4,4,0) & (4,4,0,0,0) & (4,4,0,0) \\
			\overline{5} & 1 & (1,1,1,0,0) & (1,1,1,0,0) & (1,1,1,0) & (1,1,1,1,0) & (1,1,1,0) \\
			\overline{6} & 2 & (2,2,2,0,0) & (2,2,2,0,0) & (2,2,0,0) & (2,2,2,2,0) & (2,2,2,0) \\
			\overline{7} & 3 & (3,3,3,0,0) & (3,3,3,0,0) & (3,3,0,0) & (3,3,3,0,0) & (3,3,3,0) \\
			\overline{8} & 3 & (3,3,3,0,0) & (3,3,3,0,0) & (3,3,3,0) & (3,3,3,0,0) & (3,3,0,0) \\
			\hline
		\end{array}
	\end{displaymath}
\normalsize
\end{example}

\section{The ideal $I_g$}

Once we decompose $f$ into idempotent 2-cocycles having a unique class of non-trivial annihilators we can proceed further based on the different generators of those non-trivial annihilators. Let $B^*_f$ be the set of generators of the elements of $G^*$ with respect to $f$ which are non-trivial annihilators of $f$ (for more details on the set $B_f$ see \cite{LaTh1}, Section 6). The elements of $B^*_f$ are maximal inside $\Gamma_f$ with respect to inclusion of ordered sets. For $g\in \Gamma_f$, we denote $I_g=\sum_{\sigma\in g} I_{\sigma}$.

\begin{proposition}{\label{propdecompidealgen}}
Let $f\in E^2(G,L;H)$. Then $f=\bigvee_{\gamma\in B^*_f} f_{\{J,I_{\gamma},0\}}$.
\end{proposition}

\begin{proof}
We set $f_{\gamma}=f_{\{J,I_{\gamma},0\}}$, $B^*=B^*_f$ and $H=H(f)=H(f_{\gamma})$. Since $I_{\gamma}\lhd A_f$, we know that $f_{\gamma}\leq f$, for every $\gamma\in B^*$ and so $\bigvee_{\gamma\in B^*} f_{\gamma}\leq f$. For the opposite direction, let $\sigma,\tau\in G$ such that $f(\sigma,\tau)=1$. If $f=f_0$, then the equality is immediate. So suppose that $f\neq f_0$ (i.e. there exists a generator with at least two elements). If $\sigma\in H$, then $f_{\gamma}(\sigma,\tau)=1$, for every $\gamma\in B^*$, and similarly if $\tau\in H$. If $\sigma\tau\in H$, then $\sigma,\tau\in H$ and so again $f_{\gamma}(\sigma,\tau)=1$. So let $\sigma,\tau,\sigma\tau\in G^*$. Moreover let $g_{\sigma},g_{\tau}\in\Gamma_f$, where $g_{\sigma}=(\sigma_1,\ldots,\sigma_a)$, $g_{\tau}=(\tau_1,\ldots,\tau_b)$, $a,b\geq 1$. Since $f(\sigma,\tau)=1$, we have $g_{\sigma}g_{\tau}=g_{\sigma\tau}\in\Gamma_f$ (\cite{LaTh1}, Remark 6.2). We extend $g_{\sigma\tau}$ to an element of $B^*$, say $\gamma=g_1g_{\sigma\tau}g_2$ (or $\gamma=g_1g_{\sigma\tau}$ or $\gamma=g_{\sigma\tau}g_2$ or  $\gamma=g_{\sigma\tau}$ if already $g_{\sigma\tau}\in B^*$). We note that $x_{\sigma}\in I_{\sigma_1}\subseteq I_{\gamma}$ and that $x_{\tau}\in I_{\tau_1}\subseteq I_{\gamma}$. Then $x_{\sigma\tau}=x_{\sigma}x_{\tau}\in I_{\gamma}$ and so by definition $f_{\gamma}(\sigma,\tau)=1$ which proves that $f\leq \bigvee_{\gamma\in B^*} f_{\gamma}$.
\end{proof}

\begin{example}{\label{exdecompIgamma}}
Consider the idempotent 2-cocycle $f_{I_3}$ of Example \ref{exdecomp} with $\overline{4}$ as unique non-trivial annihilator. In Example \ref{exdecompr} we saw that $f_{I_3}=f_{r_3}$. It holds that $f_{r_3}=f_{r'}$ for $r'=\{0,9,18,27,36,9,17,24,27\}$. The generators of the elements of $G^*$ with respect to $f_{r'}$ are
	\begin{eqnarray}
		\Gamma_{f_{r'}} &=& \{\{(\overline{1})\},\{(\overline{1},\overline{1})\},\{(\overline{1},\overline{1},\overline{1})\},\{(\overline{5},\overline{8}),(\overline{8},\overline{5}),(\overline{1},\overline{1},\overline{1},\overline{1})\},\{(\overline{5})\},\nonumber\\
		&& \{(\overline{6})\},\{(\overline{7})\},\{(\overline{8})\}\}.\nonumber
	\end{eqnarray}
For the three generators of $x_{\overline{4}}$, $\gamma_1=(\overline{5},\overline{8})$, $\gamma_2=(\overline{8},\overline{5})$ and $\gamma_3=(\overline{1},\overline{1},\overline{1},\overline{1})$ we have $P_1=I_{\gamma_{\overline{1}}}=I_{\gamma_{\overline{2}}}=I_{\overline{5}}+I_{\overline{8}}=Lx_{\overline{4}}+Lx_{\overline{5}}+Lx_{\overline{8}}$ and $P_2=I_{\gamma_{\overline{3}}}=Lx_{\overline{1}}+Lx_{\overline{2}}+Lx_{\overline{3}}+Lx_{\overline{4}}$. Since $(f_{r'})_{\{J,P_i,0\}}=f_{r'_{\{J,P_i,0\}}}$, $i=1,2$ (Proposition \ref{propeqcaridealr}), from Proposition \ref{propdecompidealgen} we have $f_{r'}=f_{r'_{\{J,P_1,0\}}}\vee f_{r'_{\{J,P_2,0\}}}$ with respective graphs

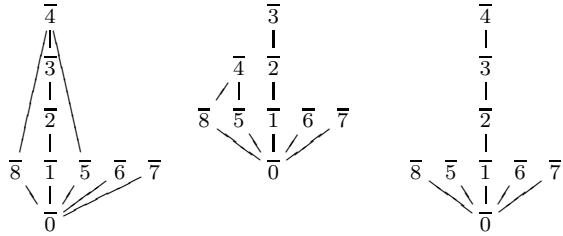
\[\scriptsize
\xymatrixrowsep{0.09in}
\xymatrixcolsep{0.04in}
\xymatrix{ & 	\overline{4}	& \\
		& \overline{3}\ar@{-}[u]& \\
				& \overline{2}\ar@{-}[u] & \\
				\overline{8}\ar@{-}[uuur]&\overline{1}\ar@{-}[u]& \overline{5}\ar@{-}[uuul] & \overline{6} & \overline{7}\\
					& \overline{0}\ar@{-}[ul]\ar@{-}[u]\ar@{-}[ur]\ar@{-}[urr]\ar@{-}[urrr] & }\quad
					\xymatrix{ & & \overline{3}& \\
				& \overline{4}& \overline{2}\ar@{-}[u]&\\
				\overline{8}\ar@{-}[ur]&\overline{5}\ar@{-}[u]&\overline{1}\ar@{-}[u]& \overline{6}& \overline{7}& \\
					& & \overline{0}\ar@{-}[ull]\ar@{-}[ul]\ar@{-}[u]\ar@{-}[ur] \ar@{-}[urr]& }\quad
					\xymatrix{  & & 	\overline{4}	& \\
					& & \overline{3}\ar@{-}[u]& \\
				& & \overline{2}\ar@{-}[u]& \\
				\overline{8} & \overline{5}& \overline{1}\ar@{-}[u]& \overline{6} &\overline{7} \\
					& & \overline{0}\ar@{-}[ull]\ar@{-}[ul]\ar@{-}[u]\ar@{-}[ur]\ar@{-}[urr]  }
\]

\noindent and generators
\small
\begin{eqnarray}
		\Gamma_{f_{r'_{\{J,P_1,0\}}}}&=& \{\{(\overline{1})\}, \{(\overline{1}, \overline{1})\}, \{(\overline{1},\overline{1},\overline{1})\}, \{(\overline{5}, \overline{8}), (\overline{8},\overline{5})\}, \{(\overline{5})\},\{(\overline{6})\},\{(\overline{7})\},\{(\overline{8})\}\}.\nonumber\\
		\Gamma_{f_{r'_{\{J,P_2,0\}}}}&=& \{\{(\overline{1})\}, \{(\overline{1}, \overline{1})\}, \{(\overline{1},\overline{1},\overline{1})\}, \{(\overline{1}, \overline{1},\overline{1},\overline{1})\}, \{(\overline{5})\},\{(\overline{6})\},\{(\overline{7})\},\{(\overline{8})\}\}.\nonumber\quad\square
  	\end{eqnarray}
\end{example}

\normalsize

Suppose that we are given an $f\in E^2(G,L)$ and we want to find some $r\in Sl(G)$ such that $f=f_r$. As a first step we should apply Proposition 6.10 of \cite{LaTh1} to every $\sigma\in G^*$ such that $\sigma$ is a non-trivial annihilator of $f$ (there is no point in applying it to trivial annihilators). For $f$ of Example \ref{exdecomp} we would have $r(\overline{5})r(\overline{8})=r(\overline{8})r(\overline{5})=r(\overline{1})^4$ and $r(\overline{1})r(\overline{1})r(\overline{5})=r(\overline{1})r(\overline{5})r(\overline{1})
=r(\overline{5})r(\overline{1})r(\overline{1})$. But if $f$ has a unique class of non-trivial annihilators, then we are left with a single equation. For $f_{I_3}$ of the same example, we would only have the first from the two mentioned above. Simpler yet, for the idempotent 2-cocycle $f_{\{J,P_1,0\}}$ of Example \ref{exdecompIgamma} we only have the equation $r(\overline{5})r(\overline{8})=r(\overline{8})r(\overline{5})$ and for $f_{\{J,P_2,0\}}$ none.

\begin{xrem}
Suppose that $f$ has a unique class of non-trivial annihilators, say $[\rho]$. If $I_{\gamma}$ is a constant for every generator $\gamma$ of $\rho$, then the equality of Proposition \ref{propdecompidealgen} is trivial, since every term $f_{\{J,I_{\gamma},0\}}$ of the second part equals $f$, as in the following example.
\end{xrem}

For $g_1,g_2\in\Gamma_f$, the relation $g_1\leq g_2$ if and only if $g_1$ is an ordered part of $g_2$, is a partial ordering with least element the empty word $()$. We call the Hasse diagram with regard to this ordering, the graph of generators of $f$.

\begin{example}
Let $D_3=\{a,b:a^3=b^2=e, bab=a^{-1}\}=\{e,a,a^2,b,ab,a^2b\}$ the dihedral group of order $6$. Let $f\in E^2(G,L)$ defined by the table

	\begin{tabular}{|c| cccccc|}
		\hline
		& $e$ & $a$ & $a^2$ & $b$ & $ab$ & $a^2b$ \\
		\hline
		$e$  & 1 & 1 & 1 & 1 & 1 & 1\\
		$a$  & 1 & 1 & 0 & 1 & 0 & 0\\ 
		$a^2$& 1 & 0 & 0 & 0 & 0 & 0\\
		$b$  & 1 & 1 & 1 & 0 & 0 & 0\\
		$ab$ & 1 & 0 & 0 & 0 & 0 & 0\\
		$a^2b$& 1 & 1 & 0 & 0 & 0 & 0\\
		\hline
\end{tabular}

\noindent with generators $\Gamma_{a}=\{(a)\}$, $\Gamma_{b}=\{(b)\}$, $\Gamma_{a^2b}=\{(b,a)\}$, $\Gamma_{a^2}=\{(a,a)\}$, $\Gamma_{ab}=\{(a,b),(b,a,a)\}$ and graphs (of generators,left,right)

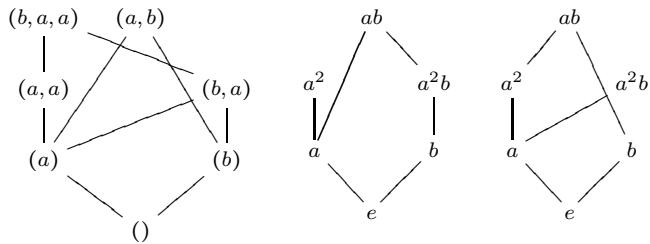
\[\scriptsize
\xymatrixrowsep{0.18in}
\xymatrixcolsep{0.1in}
\xymatrix{	 (b,a,a) & (a,b) & \\
			(a,a) \ar@{-}[u]& & (b,a) \ar@{-}[ull]  \\
				(a) \ar@{-}[u]\ar@{-}[urr] \ar@{-}[uur]& & (b) \ar@{-}[u]\ar@{-}[uul]\\
					& ()\ar@{-}[ul]\ar@{-}[ur] & &}
\xymatrix{	& ab & \\
			a^2 & & a^2b \ar@{-}[ul]  \\
				a \ar@{-}[u]\ar@{-}[uur]\ar@{-}[uur]  & & b \ar@{-}[u]\\
					& e\ar@{-}[ul]\ar@{-}[ur] & &}
\xymatrix{	& ab & \\
			a^2 \ar@{-}[ur]& & a^2b \\
				a \ar@{-}[u]\ar@{-}[urr]  & & b \ar@{-}[uul]\\
					& e\ar@{-}[ul]\ar@{-}[ur] & &}\]

\normalsize

\noindent We notice that the only non-trivial annihilator of $f$ is $ab$ with generators $\gamma_1=(a,b)$ and $\gamma_2=(b,a,a)$. Since they contain exactly the same letters ($a$ and $b$ with different multiplicities) we have $I_{\gamma_1}=I_{\gamma_2}=I_a+I_b$. No $r\in Sl(G)$ is know such that $f=f_r$. To prove that no such $r$ exists, one must prove that there does not exist a monoid $\Omega$ with the properties mentioned in the introduction such that $r(b)r(a)r(a)=r(a)r(b)$.
\end{example}

\end{document}